\documentclass[final]{siamltex}
\usepackage{amssymb,amstext,mathrsfs}
\usepackage[leqno]{amsmath}
\usepackage{cite}
\usepackage{float}
\usepackage{graphics}
\usepackage{graphicx}
\usepackage[colorlinks,linkcolor=red,anchorcolor=black,citecolor=blue]{hyperref}
\usepackage[figure,table]{hypcap}
\newtheorem{scheme}{Scheme}
\newtheorem{remark}[theorem]{Remark}

\newtheorem{thm}{Theorem}
\newcommand{\D}{\mathcal{D}}
\renewcommand{\footnote}[1]{ }

\def\F{\mathscr{F}}

\def\R{\mathbb{R}}

\def\11{\mathbf{1}}

\reversemarginpar
\def\di{\mathrm{d}}
\newcommand{\abs}[1]{\left\vert#1\right\vert}

\newcommand{\Ec}[3]{\mathbb{E}_{#1}^{#2}\left[#3\right]}

\newcommand{\HE}[1]{\mathbb{\hat{E}}^{n,x}\left[#1\right]}
\newcommand{\CE}[2]{\mathbb{E}\left[\left.#1\right\vert #2\right]}

\renewcommand{\d}{\mathrm{d}}
\newcommand{\f}{\frac}
\newcommand{\Gm}{\Gamma}
\newcommand{\rb}[2]{\raisebox{0pt}[0pt][0pt]{\raisebox{#1}{#2}}}
\newenvironment{eqsys}[1][]
{\begin{equation}#1\left\{\begin{aligned}}
{\end{aligned}\right.\end{equation}}
{\end{aligned}\end{equation}}
\newcommand{\U}{\mathcal{U}}
\newcommand{\p}{\partial}
\newcommand{\nb}{\nabla}
\newcommand{\sig}{\sigma}
\newcommand{\tp}{{\scriptscriptstyle \top}}
\renewcommand{\S}{\mathcal{S}}

 \newcommand{\bfj} {{\mathbf j}}
\renewcommand{\L}{\mathcal{L}}

\title{High order numerical schemes
for second-order FBSDEs with applications to stochastic optimal control
\thanks{This work is partially supported  by  the National
Natural Science Foundations of China under grant numbers 91130003, 11201461and 11171189, and
Natural Science Foundation of Shandong Province under grant number ZR2011AZ002.}
}

\author
{Tao Kong\thanks{School of Mathematics \& Finance Institute, Shandong University, Jinan 250100, China.
Email: vision.kt@gmail.com.}
\and Weidong Zhao\thanks{School of Mathematics \& Finance Institute, Shandong University, Jinan 250100, China.
Email: wdzhao@sdu.edu.cn.}
\and Tao Zhou\thanks{LSEC, Institute
of Computational Mathematics, Academy of Mathematics and Systems
Science, Chinese Academy of Sciences, Beijing 100190,
China. Email: tzhou@lsec.cc.ac.cn.
}}

\begin{document}

\maketitle

\begin{abstract}
This is one of our series papers on multistep schemes for solving forward backward stochastic differential equations (FBSDEs) and related problems. Here we extend (with non-trivial updates) our multistep schemes in [W. Zhao, Y. Fu and T. Zhou, SIAM J. Sci. Comput., 36 (2014), pp. A1731-A1751.] to solve the second order FBSDEs (2FBSDEs). The key feature of the multistep schemes is that the Euler method is used to discrete the forward SDE, which dramatically reduces the entire computational complexity. Moreover, it is shown that the usual quantities of interest (e.g., the solution tuple $(Y_t, Z_t, A_t, \Gamma_t)$ in the 2FBSDEs) are still of high order accuracy. Several numerical examples are given to show the effective of the proposed numerical schemes.  Applications of our numerical schemes for stochastic optimal control problems are also presented.
\end{abstract}



\pagestyle{myheadings}
\thispagestyle{plain}
\markboth{T. Kong, W. Zhao, and T. Zhou}
{High-order numerical schemes for solving second-order FBSDEs}

\section{Introduction}
This work is concerned with numerical methods for the following coupled
second order forward backward stochastic differential equations (2FBSDEs)
which is defined on the filtered probability space $(\Omega,\F, \mathbb F, P):$
\begin{equation}\label{2FBSDEs}
\left\{\begin{aligned}
X_t &= x + \int_0^t b(s,\Theta_s) \d s + \int_0^t \sigma(s,\Theta_s) \d W_s,\\
Y_t &= g(X_T) + \int_t^T f(s,\Theta_s) \d s- \int_t^T Z_s \d W_s,\\
Z_s &= Z_0 +\int_0^t A_s \d s+\int_0^t \Gamma_s \d W_s,
\end{aligned}\right. \qquad s\in [0,T],
\end{equation}
where  $\Theta_t = (X_t,Y_t,Z_t,A_t,\Gamma_t) \in \R^m\times\R\times\R^d\times \S^d $ is the unknown,
$(\Omega, \F, P)$ is the given probability space, $T>0$ is the deterministic terminal time,
 $\{W_t\}_{t\in[0,T]}$ is a $d$-dimensional Brownian Motion defined on $(\Omega, \F, P)$ with the natural filtration
$\mathbb F = \{\mathscr F_t\}_{0\le t\le T}$ and all $P$-null sets in $\F_0,$
$x\in\F_0$ is the initial condition of the
forward SDE, $\S^d$ is the set of all $d\times d$ real-valued symmetric matrices, and
\begin{align*}
b &: [0,T]\times\R^m\times\R\times\R^d\times \S^d \rightarrow\R^m, &
\sigma & : [0,T]\times\R^m\times\R\times\R^d\times \S^d \rightarrow\R^{m\times d}
\end{align*}
are referred to the drift and  diffusion coefficients of the forward SDE, respextively. While
\begin{align*}
f &: [0,T]\times\R^m\times\R\times\R^d\times \S^d\rightarrow\R, &
g&: \R^m \to \R
\end{align*}
are referred to the generator and the terminal condition of
the backward SDE, respectively. The two stochastic integrals with respect to the Brownian Motion $\{W_t\}_{t\in[0,T]}$
are of the It\^o type. A 5-tuple $(X_t,Y_t,Z_t,\Gamma_t,A_t)$  is called an $L^2$-adapted solution
of the 2FBSDEs (\ref{2FBSDEs}) if
it is $\F_t$-adapted, square integrable, and satisfies \eqref{2FBSDEs}.
Moreover, the 2FBSDE \eqref{2FBSDEs} is called \textit{decoupled} if $b$ and $\sigma$
are independent of $Y_t$, $Z_t$ $A_t$ and $\Gamma_{t}$.

The 2FBSDEs (\ref{2FBSDEs}), as an extension of the BSDEs in the linear \cite{Bismut} or nonlinear cases \cite{PP1990}, was first introduced by the cornerstone work of Cheridito, Soner, Touzi and Victoir \cite{CSTV}, with a slightly different (yet equivalent) formula, and it was further investigated by Soner, Touzi and Zhang in \cite{STZ}. The main motivation there is to give a precise connection between the 2FBSDEs and fully non-linear PDEs, in particular
the Hamilton-Jacobi-Bellman equations and the Bellman-Isaacs equations
which are widely used in stochastic control and in stochastic differential
games. Such a connection leads to interesting stochastic representation results for fully nonlinear PDEs, generalizing the original (nonlinear) Feynman-Kac representations of linear and semi-linear parabolic
PDEs, see e.g. \cite{KKPPQ,PT,Peng,MY} and references therein.

From the numerical point of view, one can adopt these connections between PDEs and FBSDEs to design the so called probabilistic numerical methods for PDEs, by solving the equivalent FBSDEs (or 2FBSDEs). While there are a lot of work dealing with numerical schemes for BSDEs \cite{BT,BZ,JF,Cubature,LGW,Z1,ZZJ,ZCP}, however, there are only a few work on numerical methods for FBSDEs \cite{DMP,MT,MPJ,MSZ,FZZ2015,TZZ2015,ZFZ2014,ZZwJ} and fully non-linear PDEs \cite{FTW,GZZ}.
Some of the above works are designed with high order accuracy that can however only be used to deal with low dimensional FBSDEs. While some of them are low order numerical methods that are suitable for solving high dimensional problems. In particular, we mention the work \cite{GZZ}, where a numerical example for a 12-dimensional coupled FBSDE is reported, and it is shown by numerical test that the numerical method converges with order 1.  Also, in \cite{FZZ2015}, multistep schemes were proposed to solve multi-dimensional FBSDEs by using the spares grid interpolation technique, and several multi-dimensional examples with dimension up to 6 are presented, and high order convergence rates up to 3 were obtained.

To the best of our knowledge, there is no related studies for high order numerical methods for 2FBSDEs.
The main purpose in this work is to extend our multistep schemes in \cite{ZFZ2014} (which is original designed for solving FBSDEs) to the use of solving 2FBSDEs. The key feature of the multistep schemes is that the Euler method is used to discrete the forward SDE, which dramatically reduces the entire computational complexity. Moreover, it is shown that the quantities of interest (e.g., the solution tuple $(Y_t,Z_t,\Gamma_t,A_t)$ in the 2FBSDEs) are still of high order accuracy. Several numerical examples are given to show the effective of the proposed numerical schemes.  With particular interests, applications of the proposed numerical schemes to stochastic optimal control problems are also presented and corresponding numerical examples are also given. It is worth to note that the multistep schemes proposed in this paper can also be used to solve a large class of fully nonlinear PDEs and we left this for our future studies.

The rest of the paper is organized as follows. In Section 2, we present some preliminaries. Then in Section 3, We shall discuss the multistep schemes for solving decoupled 2FBSDEs, and this is followed by extensions to coupled 2FBSDEs in Section 4. In Section 5, we shall present several numerical examples, and applications of our numerical methods for stochastic optimal control problems will also be presented. We finally give some concluding remarks in Section 6.

\section{Preliminaries}\label{pre}

We first introduce some notations that will be used in this paper. For
$x\in \R^d$,  $ B \in \S^d$, we set $\abs{x} := \sqrt{x_1^2 +\dots + x_d^2}$,
 $\abs{B} := \sqrt{\sum_{i,j=1}^d B_{ij}^2}$.
We denote by $C_b^k$ the set of functions $g(x):\R\to\R$ with uniformly bounded derivatives up to order $k$,
and by $C^{k_1,k_2}$ the set of functions $f(t,x):[0,T]\times\R^n\to\R^m$ with
continuous partial derivatives up to order $k_1$ w.r.t.\ $t\in\R$ and up to order $k_2$ w.r.t.\ $x\in\R^m$.

\subsection{Diffusion process and its generator}
A stochastic process $X_t$ is called a diffusion process starting at $x_0$ and at the time $t_0$
if it satisfies the following SDE
\begin{equation}\label{SDEs}
\begin{aligned}
X_t =& x_0 + \int_{t_0}^t b(s,X_s) \di s + \int_{t_0}^t \sigma(s,X_s) \di W_s, \quad t\in[t_0, T],
\end{aligned}
\end{equation}
where $b_s=b(s,X_s)$ and $\sigma_s=\sigma(s,X_s)$ are measurable functions satisfying
\begin{equation}\label{lcsde}\begin{aligned}
\abs{b(s,x)}+\abs{\sigma(s,x)} &\leq C(1+\abs{x}), \quad x\in \R^m, s\in [t_0,T],\\
\abs{b(t,x)-b(t,y)} + \abs{\sigma(t,x)-\sigma(t,y)} &\leq L\abs{x-y}, \quad x,y\in \R^m, s\in [t_0,T].
\end{aligned}\end{equation}
It is well known that under conditions \eqref{lcsde}, the SDE \eqref{SDEs}
admits a unique solution. It is worth to note that, by the Markov property of the diffusion process, we have
$\Ec{t}{x}{X_s} = \CE{X_s}{X_t=x}, \forall t\leq s$.

Let $X_t$ be the solution of \eqref{SDEs}. Then, for any given measurable function $g: [0,T]\times \R^m\rightarrow \R$,
$g(t,X_t)$ becomes a stochastic process. Moreover, we give the following definition
\begin{definition}\label{generator}
Let $X_s$ be a diffusion process in $\R^m$ that satisfies \eqref{SDEs}.
The generator $A$ of $X_s$ on $g$ is defined by
\begin{equation}\label{gnrt}
A g(t,x) = \lim_{s\downarrow t}\frac{\Ec{t}{x}{g(s,X_s)}-g(t,x)}{s-t},\quad x\in\R^n.
\end{equation}
\end{definition}
Concerning the generator $A,$ the following result holds \cite{O}:
\begin{thm}\label{thm_generator}
Let $X_s$ be the diffusion process defined by the SDE (\ref{SDEs}).
If $f\in C^{1,2}([0,T]\times\R^m)$, then we have
\begin{equation}\label{gnrt_dif}
A f(t,x) = \L f(t,x), \quad A f(t,X_t) = \L f(t,X_t),
\end{equation}
where the operator  $\L$ is defined by
\begin{equation}\label{operatorL}
\L\phi(t,x) := \phi_t(t,x) + \nb_x\phi(t,x)b(t,x) +
\frac12\mathrm{tr}\left(\sigma(t,x)\sigma^\tp(t,x)\nb_x^2\phi(t,x)\right).
\end{equation}
with $\nb_x\phi = (\p_{x_1}\phi,\dots,\p_{x_m}\phi)$,
 and $\nb_x^2\phi$ being the Hessian matrix of $\phi$ with respect to the spatial variable $x$.
\end{thm}

Note that $A f(t,X_t)\in \mathcal F_t$ is a stochastic process. Furthermore, by using together the It\^o's formula, Theorem \ref{thm_generator}
and the tower rule of conditional
expectations, we have the following theorem.


%

\begin{thm}\label{thm_identity}
Let $t_0<t$ be a fixed time, and $x_0\in \R^m$ be a fixed space point. If $f\in C^{1,2}([0,T]\times\R^m)$
and $\Ec{t_0}{x_0}{\abs{A f(t,X_t)}}<+\infty$,
we have
\begin{equation}\label{gnrt_eq}
\frac{\di\Ec{t_0}{x_0}{f(t,X_t)}}{\di t} = \Ec{t_0}{x_0}{Af(t,X_t)}, \quad t\ge t_0.
\end{equation}
Moreover, the following identity holds
\begin{equation}\label{gnrt_eq1}
\left. \frac{\di\Ec{t_0}{x_0}{f(t,X_t)}}{\di t} \right|_{t=t_0} = \left.\frac{\di\Ec{t_0}{x_0}{f(t,\bar X_t)}}{\di t}\right|_{t=t_0},
\end{equation}
where $\bar X_t$ is a diffusion process satisfying
\begin{equation}
\bar{X}_t = x_0 + \int_{t_0}^t\bar{b}_s\di s + \int_{t_0}^t\bar\sigma_s\di W_s,
\end{equation}
with
$\bar b:[0,T]\times\R^m\to\R^m, \bar\sig:[0,T]\times\R^m\to\R^{m\times d}$ being smooth functions 
satisfying
\begin{equation*}
\bar b(t_0,x_0) = b(t_0,x_0), \quad \bar \sig(t_0,x_0) = \sig(t_0,x_0).
\end{equation*}
\end{thm}
Note that by choosing different $\bar b$ and $\bar\sigma$,
the identity \eqref{gnrt_eq1} yields different ways for approximating
$\left.\di\Ec{t_0}{x}{f(t,X_t)}/\di t \right|_{t=t_0}$.
The computational complexity can be significantly reduced if
appropriate choices are made.

\subsection{Solution regularity and representation theory of 2FBSDEs}
We now consider the following decoupled 2FBSDEs
\begin{eqsys}\label{dp2FBSDEs}
\d X_t &= b(t,X_t) \d t + \sig(t,X_t) \d W_t,\\
-\d Y_t &= f(t,X_t, Y_t, Z_t,\Gm_t) \d t - Z_t  \d W_t,\\
\d Z_t &= A_t \d t + \Gm _t \d W_t\\
\end{eqsys}%
with a given terminal condition $Y_T = g(X_T).$  Under some standard assumptions (for details, please refer to \cite{CSTV, STZ}),
the following results are shown \cite{CSTV}:
\begin{theorem}
Let $u=u(t,x)$ be the solution of the following fully nonlinear PDE
 \begin{equation}\begin{aligned}\label{pde}
 \L u + f(t,x,u,\nabla_x u\sigma, \nabla_x(\nabla_x u\sigma)\sigma) =0
 \end{aligned}\end{equation}
with the terminal condition $u(T,x) = g(x)$, and let $(X_t, Y_t, Z_t, \Gm_t, A_t)$ be the solution of
the 2FBSDE (\ref{dp2FBSDEs}). Then we have
\begin{align*}
&\qquad Y_t = u(t,X_t), \quad Z_t =(\nabla_x u \sigma)(t,X_t), \\
& \Gamma_t =  (\nabla_x (\nabla_x u \sigma)\sigma)(t,X_t), \quad A_t =  (\L (\nabla_x u \sigma))(t,X_t).
\end{align*}
where the associate operator $\L$ is defined by \eqref{operatorL}.
\end{theorem}

The above theorem provides a stochastic representation for solutions of fully nonlinear parabolic PDEs, generalizing the pioneer work on Feynman-Kac representations of linear and semi-linear parabolic
PDEs \cite{KKPPQ,PT,Peng,MY}.

\subsection{Derivative approximation}\label{app_diff}
We now recall some basic results for numerical approximation of derivatives, and these results will play an
important role for designing our high order numerical methods for 2FBSDEs.

Let $u(t)\in C_b^{k+1}$ and $t_i\in \mathbb R$ $(0\le i\le k)$ satisfying $t_0<t_1<\cdots<t_k$,
where $k$ is a positive integer. Let $\Delta t_{0,i}=t_i-t_0,i=0,1,\ldots,k$.
Then by Taylor's expansion, for each $t_i,i = 0,1,\ldots,k$,
we have
\begin{flalign*}
\begin{split}
u(t_i) &= \sum_{j=0}^k\frac{(\Delta t_{0,i})^j}{j!}\frac{d^j u}{dt^j}(t_0) + O\left(\Delta t_{0,i}\right)^{k+1}.
\end{split}
\end{flalign*}
Then, we can deduce
\begin{flalign*}
\begin{split}
\sum_{i=0}^k \alpha_{k,i} u(t_i) &= \sum_{j=0}^k\frac{\sum\limits_{i=0}^k\alpha_{k,i}(\Delta t_{0,i})^j}{j!}\frac{d^j u}{dt^j}(t_0)
+O\left(\sum_{i=0}^k\alpha_{k,i} (\Delta t_{0,i})^{k+1}\right),
\end{split}
\end{flalign*}
where $\alpha_{k,i}$, $i=0,1,\dots,k$, are real nubmers.
By choosing $\alpha_{k,i}~(i=0,1,\dots,k)$ as
\begin{equation}\label{alpha}
\frac{\sum\limits_{i=0}^k\alpha_{k,i} (\Delta t_{0,i})^j}{j!}=
 \begin{cases}
 1, & j=1\\
 0, & j\neq 1
 \end{cases}.
\end{equation}
One obtains a high order approximation
\begin{equation}\label{dire_appro}
\frac{d u}{dt}(t_0) = \sum_{i=0}^k\alpha_{k,i} u(t_i) + R_D,
\end{equation}
where $R_D = O(\sum_{i=0}^k\alpha_{k,i} (\Delta t_{0,i})^{k+1})$.
In particular, when $\Delta t_{0,i} = i\Delta t$, we get from \eqref{alpha} the following linear system for $\alpha_{k,i}\Delta t$,
\begin{equation}\label{alpha_equ}
 \sum_{i=1}^k i^j[\alpha_{k,i}\Delta t]  = 
 \begin{cases}
 1, & j=1\\
 0, & j\neq 1
 \end{cases}.
 \end{equation}
Note that the above system can be solved easily. We list $\alpha_{k,i}\Delta t$ $(i=0,1,\dots,k)$
of the system \eqref{alpha_equ} for $k=1,2,\ldots,6$ in the following table.

\begin{table}[H]\label{alpha_value}
\renewcommand{\arraystretch}{1.25}
\setlength{\belowcaptionskip}{0pt}
\caption{}
\footnotesize
\begin{center}
\begin{tabular}{c|c|c|c|c|c|c|c}
\hline
$\alpha_{k,i}\Delta t$ & $i=0$ & $i=1$ & $i=2$ & $i=3$ & $i=4$ & $i=5$ & $i=6$\\
\hline
$k=1$ & $-1$              & $1$ &                 &                &                 &           & \\
\hline
$k=2$ & $-\frac32$        & $2$ & $-\frac12$      &                &                 &           & \\
\hline
$k=3$ & $-\frac{11}{6}$   & $3$ & $-\frac32$      & $\frac13$      &                 &           & \\
\hline
$k=4$ & $-\frac{25}{12}$  & $4$ & $-3$            & $\frac43$      & $-\frac14$      &           & \\
\hline
$k=5$ & $-\frac{137}{60}$ & $5$ & $-5$            & $\frac{10}{3}$ & $-\frac54$      & $\frac15$ & \\
\hline
$k=6$ & $-\frac{49}{20}$  & $6$ & $-\frac{15}{2}$ & $\frac{20}{3}$ & $-\frac{15}{4}$ & $\frac65$ & $-\frac16$\\
\hline
\end{tabular}\end{center}
\end{table}
The above approximation schemes are very popular in numerical ODE community, and it is known that
such an approximation scheme is unstable for $k\geq 7,$ and this is why we have only listed the values of $\alpha_{k,i}\Delta t$ for $1\le k\le 6$ in Table \ref{alpha_value}. For more details, one can refer to \cite{ZFZ2014}.

\section{Multistep schemes for decoupled 2FBSDEs}
In this section, we first consider the multistep schemes for solving decoupled 2FBSDEs \eqref{dp2FBSDEs}, i.e., the forward SDE is independent of $(Y_t, Z_t, A_t, \Gamma_t).$  To begin, let $N$ be a positive integer. For the time interval $ [t_0,T],$ we introduce a regular time partition:
$$
t_0 < t_1 <\dots<t_N = T.
$$
We set $\Delta t_{t_n,k} = t_{n+k} - t_n$ for $n,\in \{1,2\dots,N\}$ and $k \in \mathbb{N} $ satisfying $n+k \le N$.
For $t \ge t_n$, we denote  $\Delta W_{t_n,k} = W_{t_{n+k}} - W_{t_n},$ $\Delta t_{t_n,t} =t- t_n$ and $\Delta W_{t_n,t} = W_t - W_{t_n}$,

\subsection{Four reference ODEs}
Let $\Theta_{t}=(X_t,Y_t,Z_t,A_t,\Gamma_t)$ be the solution of the decoupled FBSDEs \eqref{dp2FBSDEs} with terminal condition $g(X_T)$. We set $\Ec{t_n}{x}{\cdot} = \CE{\cdot}{\F_{t_n}^{t_n,x}}$.
By taking conditional expectation $\Ec{t_n}{x}{\cdot}$ on both sides of the last two equations in \eqref{dp2FBSDEs},
we obtain
\begin{align}
\Ec{t_n}{x}{Y_t} =  &\Ec{t_n}{x}{g(X_T)} + \int_t^T \Ec{t_n}{x}{f(s,\Theta_s)}\di s, &t \in [t_n, T],\label{exp_y} \\
\Ec{t_n}{x}{Z_t} =  & \Ec{t_n}{x}{Z_{t_n}}  +
\int_{t_n}^t  \Ec{t_n}{x}{A_s} ds, &t \in [t_n, T],\label{exp_a}
\end{align}
By taking derivative with respect to $t$ in \eqref{exp_y}-\eqref{exp_a}, one gets the following two reference ODEs:
\begin{align}
\frac{\di \Ec{t_n}{x}{Y_t}}{\di t} & = -\Ec{t_n}{x}{f(t,\Theta_t)},&t\in[t_n,T],\label{dif_y}\\
\frac{\di\Ec{t_n}{x}{Z_t}}{\d t} & =   \Ec{t_n}{x}{A_t},
& t\in[t_n,T]\label{dif_a}
\end{align}
Under sufficient regularity assumptions on the given data, the integrand $\Ec{t_n}{x}{f(s,\Theta_s)}$ and $\Ec{t_n}{x}{A_t}$ are continuous at $s=t$. Note that we also have
\begin{equation*}
Y_{t_n} = Y_t + \int_{t_n}^t f(s,\Theta_s) \di s - \int_{t_n}^t Z_s \di W_s,
\quad t \in [t_n,T]
\end{equation*}
By multiplying both sides of the above equation
by $\Delta W_{t_n,t}^\tp$ (where "$^\tp$" stands for the transposition),
and taking the conditional expectation $\Ec{t_n}{x}{\cdot}$ on both sides of the derived equation,
we obtain, for $t \in [t_n,T]$
\begin{flalign}\label{exp_z}
\begin{split}
0 =~& \Ec{t_n}{x}{Y_{t}\Delta W_{t_n,t}^\tp}
+\int_{t_n}^{t}\Ec{t_n}{x}{f(s,\Theta_s)\Delta W_{t_n,s}^\tp}\di s
- \int_{t_n}^{t}\Ec{t_n}{x}{Z_s}\di s,
\end{split}
\end{flalign}
Similarly, for the last equation in \eqref{dp2FBSDEs}, we have for $t \in [t_n,T]$
\begin{equation}\label{exp_gam}
0 = \Ec{t_n}{x}{Z_t^\tp\Delta W_{t_n,t}^\tp} - \int_{t_n}^{t}\Ec{t_n}{x}{A_s^\tp \Delta W_{t_n,t}^\tp}\d s
-\int _{t_n}^{t}\Ec{t_n}{x}{\Gm_s} \d s.
\end{equation}
Now, by taking derivative with respect to $t\in[t_n,T)$ on both sides,
one gets the following reference ODEs:
\begin{align}
\frac{\di \Ec{t_n}{x}{Y_t\Delta W_{t_n,t}^\tp}}{\di t} &= -\Ec{t_n}{x}{f(t,\Theta_t)\Delta W_{t_n,t}^\tp}
+ \Ec{t_n}{x}{Z_t},& t\in[t_n,T]\label{dif_z}\\
\frac{\di \Ec{t_n}{x}{Z_t^\tp\Delta W_{t_n,t}^\tp}}{\di t} &=
  \Ec{t_n}{x}{A_t^\tp\Delta W_{t_n,t}^\tp}
  + \Ec{t_n}{x}{\Gamma_t}, & t\in[t_n,T]\label{dif_gam}
\end{align}
The equations \eqref{dif_y}, \eqref{dif_a}, \eqref{dif_z}, and \eqref{dif_gam}
are reference ODEs for the decoupled 2FBSDEs \eqref{dp2FBSDEs}. Our multistep numerical schemes
will be constructed by approximating the associate derivatives and the conditional expectations
in these ODEs.

\subsection{The semi-discrete scheme}
Motivated by Theorem \ref{thm_identity}, we choose smooth functions
$\bar b(t,x)$ and $\bar\sigma(t,x)$ for $t\in [t_n,T]$ and $x\in \mathbb R^m$
with the constraints $\bar b(t_n,x)=b(t_n,x)$ and
$\bar \sigma(t_n,x)=\sigma(t_n,x)$.
Let $\bar X_t^{t_n,x}$ be the diffusion process defined by
\begin{equation}\label{xbar}\begin{aligned}
\bar{X}_t^{t_n,x} & = x + \int_{t_n}^t \bar b(s,\bar X_s^{t_n,x})\di s + \int_{t_n}^t \bar\sigma(s,\bar X_s^{t_n,x})\di W_s.
\end{aligned}\end{equation}

Note that $Y_t,Z_t,\Gm_t,A_t$ are all functions of $(t,X_t).$ Now, let $(\bar Y_t^{t_n,x},\bar Z_t^{t_n,x})$
be the value of function $(Y_t, Z_t)$ at the time-space point $(t,\bar{X}_t^{t_n,x})$. ,
by Theorem \ref{thm_identity}, we have
\begin{align*}
\left. \frac{\di\Ec{t_n}{x}{Y_t}}{\di t} \right|_{t=t_n}
&= \left.\frac{\di\Ec{t_n}{x}{\bar Y_t^{t_n,x}}}{\di t}\right|_{t=t_n},&
\left. \frac{\di \Ec{t_n}{x}{Y_t\Delta W_{t_n,t}^\tp}}{\di t}\right|_{t=t_n}
&= \left.\frac{\di \Ec{t_n}{x}{\bar Y_t^{t_n,x}\Delta W_{t_n,t}^\tp}}{\di t}\right|_{t=t_n},\\
\left. \frac{\di\Ec{t_n}{x}{Z_t}}{\di t} \right|_{t=t_n}
&= \left.\frac{\di\Ec{t_n}{x}{\bar Z_t^{t_n,x}}}{\di t}\right|_{t=t_n},&
\left. \frac{\di \Ec{t_n}{x}{Z_t^\tp\Delta W_{t_n,t}^\tp}}{\di t}\right|_{t=t_n}
&= \left.\frac{\di \Ec{t_n}{x}{(\bar Z_t^{t_n,x})^\tp\Delta W_{t_n,t}^\tp}}{\di t}\right|_{t=t_n}.
\end{align*}
Now applying
\eqref{dire_appro} to
$\left.\frac{\di\Ec{t_n}{x}{\bar Y_t^{t_n,x}}}{\di t}\right|_{t=t_n}$ and
$\left.\frac{\di \Ec{t_n}{x}{\bar Y_t^{t_n,x}\Delta W_{t_n,t}^\tp}}{\di t}\right|_{t=t_n}$,
$\left.\frac{\di\Ec{t_n}{x}{\bar Z_t^{t_n,x}}}{\di t}\right|_{t=t_n}$ and
$\left.\frac{\di \Ec{t_n}{x}{(\bar Z_t^{t_n,x})^\tp\Delta W_{t_n,t}^\tp}}{\di t}\right|_{t=t_n}$,
we deduce
\begin{equation}\label{app}
\left\{
\begin{aligned}
\left.\frac{\di \Ec{t_n}{x}{Y_t}}{\di t}\right\vert_{t=t_n}
&= \sum_{i=0}^k\alpha_{k,i}\Ec{t_n}{x}{\bar Y_{t_{n+i}}^{t_n,x}} + \bar{R}_{y,n}^k,\\
\left.\frac{\di \Ec{t_n}{x}{Y_t\Delta W_{t_n,t}^\tp}}{\di t}\right\vert_{t=t_n}
&= \sum_{i=1}^k \alpha_{k,i} \Ec{t_n}{x}{\bar Y_{t_{n+i}}^{t_n,x}\Delta W_{n,i}^\tp}
+\bar{R}_{z,n}^k, \\
\left.\frac{\di \Ec{t_n}{x}{Z_t}}{\di t}\right\vert_{t=t_n}
&= \sum_{i=0}^k\alpha_{k,i}\Ec{t_n}{x}{\bar Z_{t_{n+i}}^{t_n,x}} + \bar{R}_{A,n}^k,\\
\left.\frac{\di \Ec{t_n}{x}{Z_t^\tp\Delta W_{t_n,t}^\tp}}{\di t}\right\vert_{t=t_n}
&= \sum_{i=1}^k \alpha_{k,i} \Ec{t_n}{x}{(\bar Z_{t_{n+i}}^{t_n,x})^\tp\Delta W_{n,i}^\tp}
+\bar{R}_{\Gamma,n}^k,
\end{aligned}\right.
\end{equation}
Here $\alpha_{k,i}$ are given by \eqref{alpha},
and $\bar{R}_{y,n}^k$, $\bar{R}_{z,n}^k$,
$\bar{R}_{A,n}^k$ and $\bar{R}_{\Gamma,n}^k$ are the associate truncation errors, i.e.
\begin{align*}
\bar{R}_{y,n}^k &=\left.\frac{\di \Ec{t_n}{x}{Y_t}}{\di t}\right\vert_{t=t_n}
-\sum_{i=0}^k\alpha_{k,i}\Ec{t_n}{x}{\bar Y_{t_{n+i}}^{t_n,x}}, \\
\bar{R}_{z,n}^k &= \left.\frac{\di \Ec{t_n}{x}{Y_t\Delta W_{t_n,t}^\tp}}{\di t}\right\vert_{t=t_n}
- \sum_{i=1}^k \alpha_{k,i} \Ec{t_n}{x}{\bar Y_{t_{n+i}}^{t_n,x}\Delta W_{n,i}^\tp}, \\
\bar{R}_{A,n}^k &=\left.\frac{\di \Ec{t_n}{x}{Z_t}}{\di t}\right\vert_{t=t_n}
-\sum_{i=0}^k\alpha_{k,i}\Ec{t_n}{x}{\bar Z_{t_{n+i}}^{t_n,x}}, \\
\bar{R}_{\Gamma,n}^k &= \left.\frac{\di \Ec{t_n}{x}{Z_t^\tp\Delta W_{t_n,t}^\tp}}{\di t}\right\vert_{t=t_n}
- \sum_{i=1}^k \alpha_{k,i} \Ec{t_n}{x}{(\bar Z_{t_{n+i}}^{t_n,x})^\tp\Delta W_{n,i}^\tp}.
\end{align*}
By inserting \eqref{app} into \eqref{dif_y}, \eqref{dif_a},\eqref{dif_z} and \eqref{dif_gam},
respectively, we deduce
\begin{equation}\label{ref_all}
 \left\{
\begin{aligned}
&\sum_{i=0}^k\alpha_{k,i}\Ec{t_n}{x}{\bar Y_{t_{n+i}}} = -f(t_n,x,Y_{t_n},Z_{t_n}) + R_{y,n}^{k},\\
&\sum_{i=1}^k\alpha_{k,i}\Ec{t_n}{x}{\bar Y_{t_{n+i}}\Delta W_{n,i}^\tp} = Z_{t_n} + R_{z,n}^k\\
&\sum_{i=0}^k\alpha_{k,i}\Ec{t_n}{x}{\bar Z_{t_{n+i}}} = A_{t_n} + R_{A,n}^{k},\\
&\sum_{i=1}^k\alpha_{k,i}\Ec{t_n}{x}{\bar Z_{t_{n+i}}^\tp\Delta W_{n,i}^\tp} = \Gamma_{t_n} + R_{\Gamma,n}^k,
\end{aligned}\right.
\end{equation}
where $R_{y,n}^k = -\bar R_{y,n}^k, R_{z,n}^k= -\bar R_{z,n}^k, R_{\Gamma,n}^k = -\bar R_{\Gamma,n}^k$ and $R_{A,n}^k = -\bar R_{A,n}^k$.

Let $Y^n$, $Z^n$, $A^n$ and $\Gamma^n$
be four random variables that represent the approximate
values of the solutions $Y_t$, $Z_t$, $A_t$ and $\Gamma_t$
of the 2FBSDE in \eqref{dp2FBSDEs} at time $t_n,$ respectively.
By removing the truncation error terms $R_{y,n}^k$, $R_{z,n}^k$,
$R_{A,n}^k$ and $R_{\Gamma,n}^k$
from \eqref{ref_all}, we obtain the following semi-discrete numerical scheme
for solving the 2FBSDEs \eqref{dp2FBSDEs}:
\begin{scheme}\label{semi_sch}
Assume that random variables $Y^{N-i}$ and $Z^{N-i}$, $i = 0,1,\ldots,k-1$ are known.
For $n=N-k,\ldots,0$, with $X_t^{t_n,X^n}$ being the solution of \eqref{xbar}, solve  $Y^{n}=Y^n(X^n)$,\, $Z^{n}=Z^n(X^n)$, \,$A^{n}=A^n(X^n)$\, and $\Gamma^{n}=\Gamma^n(X^n)$  by
\begin{align}
Z^n &= \sum_{j=1}^k\alpha_{k,j}\Ec{t_n}{X^n}{\bar Y^{n+j}\Delta W_{n,j}^\tp},\label{semi_sch_z}\\
\Gamma^n &= \sum_{j=1}^k\alpha_{k,j}\Ec{t_n}{X^n}{(\bar Z^{n+j})^\tp\Delta W_{n,j}^\tp},\label{semi_sch_Ga}\\
A^n &= \sum_{j=0}^k\alpha_{k,j}\Ec{t_n}{X^n}{\bar Z^{n+j}},\label{semi_sch_A},\\
-\alpha_{k,0}Y^n &= \sum_{j=1}^k\alpha_{k,j}\Ec{t_n}{X^n}{\bar Y^{n+j}}
+f(t_n,X^n,Y^n,Z^n,\Gamma^n),\label{semi_sch_y}
\end{align}
where $\bar Y^{n+j}$ and $\bar Z^{n+j}$ are the values of
$Y^{n+j}$ and $Z^{n+j}$ at the space points $X_{t_{n+j}}^{t_n,X^n}$.
\end{scheme}

In addition, we choose $\bar b$ and $\bar \sigma$ in \eqref{xbar} as
\begin{equation}
\bar b(s,X_s^{t_n,x}) = b(t_n,x), \quad
\bar \sigma(s,X_s^{t_n,x}) = \sigma(t_n,x),
\quad \forall s\in[t_n, T].
\end{equation}
Then, \eqref{xbar} yields the Euler discretization scheme for the forward SDE.
In this case, Scheme \ref{semi_sch} becomes
\begin{scheme}\label{semi_sch_Eul}
Assume that random variables $Y^{N-i}$ and $Z^{N-i}$,
 $i = 0,1,\ldots,k-1$, are known.
For $n=N-k,\ldots,0$, solve
the random variables $X^{n,j},$ $(j=1,2,\dots,k)$, $Y^n$,\, $Z^n$, \,$A^n$ and $\Gamma^n$ by
\begin{align}
X^{n,j} &= X^n + b(t_n, X^n)\Delta t_{n,j} + \sigma(t_n, X^n)\Delta W_{n,j},\quad j = 1,\ldots,k,\label{semi_sch_Eul_x}\\
Z^n &= \sum_{j=1}^k\alpha_{k,j}\Ec{t_n}{X^n}{\bar Y^{n+j}\Delta W_{n,j}^\tp},\label{semi_Eul_z}\\
A^n &= \sum_{j=0}^k\alpha_{k,j}\Ec{t_n}{X^n}{\bar Z^{n+j}},\label{semi_Eul_A}\\
\Gamma^n &= \sum_{j=1}^k\alpha_{k,j}\Ec{t_n}{X^n}{(\bar Z^{n+j})^\tp\Delta W_{n,j}^\tp},\label{semi_Eul_Ga}\\
-\alpha_{k,0}Y^n &= \sum_{j=1}^k\alpha_{k,j}\Ec{t_n}{X^n}{\bar Y^{n+j}}
+f(t_n,X^n,Y^n,Z^n,\Gamma^n),\label{semi_Eul_y}
\end{align}
where $\bar Y^{n+j}$ and $\bar Z^{n+j}$
stand for the values of
$Y^{n+j}$ and $Z^{n+j}$ at the space points $X^{n,j}$, respectively.
\end{scheme}

If the functions
$\Ec{t_n}{x}{Y_t}$, $\Ec{t_n}{x}{Y_t\Delta W_{t_n,t}^\tp}$,
$\Ec{t_n}{x}{Z_t}$, $\Ec{t_n}{x}{Z_t^\tp\Delta W_{t_n,t}^\tp}$ and their derivatives (with respect to $t$) up to order $k+1$
are bounded, we have the estimates
\begin{equation}\label{trun_error_de}
\bar R_{y,n}^k=O\left(\Delta t\right)^k,\quad \bar R_{z,n}^k = O\left(\Delta t\right)^k,
\quad \bar R_{A,n}^k=O\left(\Delta t\right)^k,\quad \bar R_{\Gamma,n}^k = O\left(\Delta t\right)^k.
\end{equation}

\subsection{The fully discrete scheme}\label{Full_sch}

We now propose our fully discrete schemes. To this end,
we first  introduce  the time-space  partition.  Given a time partition $\mathcal{T}:=\{t_0,t_1,\dots,t_N\}$, we introduce a series of space partitions $\mathcal{D}_{\mathbf{h}} := \{\D_{h_n}^n\}_{n=0,1, \dots,N}$ for space $\R^m$,  with each $\D_{h_n}^n$ corresponds to the time level $t_n\in \mathcal{T}$ and  $\mathbf{h}=(h_0,h_1,\dots,h_N)$ being the partition density of each $\D_{h_n}^n$. We use $\D_h^n$ for simplicity if the context is clear. In addition, If  the space partitions w.r.t.\ all the  time levels  are the same, then $\D_h$ is used to represent the unified space partition.

For the space partition  $\D_{h_n}^n=\{x_j\in\R^m|j\in I, I\subset \mathbb{N}\}$, $x_j$ is called the grid points and $I$ is called the index set of the grid.  The partition density $h_n$ of  $\D_{h_n}^n$ is defined by
$h_n = \max\limits_{x_i,x_j\in\D_{h_n}^n } d(x_i,x_j)$, with $x_i, x_j$ adjacent  in $\D_{h_n}^n$, where $d(\cdot,\cdot)$ is the distance between two points in $\R^m$.
Given $x\in\R^m$, a neighbor gird  $U_{n,x}^{\D}$ of $x$, is a finite-element  subset of  $\D_{h_n}^n$ satisfying
 $\min\limits_{x_i\in U_{n,x}^{\D}} d(x,x_i) < \min\limits_{x_i\in \D_{h_n}^n \backslash U_{n,x}^{\D}} d(x,x_i)$.

With the above spacial partition, we seek to solving
$Y^n$, $Z^n$, $A^n$ and $\Gamma^n$
at grid point $x\in \D_h^n.$
More precisely, for each $x\in \D_{h}^n,n=N-k,\ldots,0$, we aim at solving $Y^n$, $Z^n$, $\Gamma^n$ and $A^n$ by
\begin{equation}\label{semi_Eul}
\left\{
\begin{aligned}
Z^n &= \sum_{j=1}^k\alpha_{k,j}\Ec{t_n}{x}{\bar Y^{n+j}\Delta W_{n,j}^\tp},\\
\Gamma^n &= \sum_{j=1}^k\alpha_{k,j}\Ec{t_n}{x}{(\bar Z^{n+j})^\tp\Delta W_{n,j}^\tp},\\
A^n &= \sum_{j=0}^k\alpha_{k,j}\Ec{t_n}{x}{\bar Z^{n+j}},\\
-\alpha_{k,0}Y^n &= \sum_{j=1}^k\alpha_{k,j}\Ec{t_n}{x}{\bar Y^{n+j}}
+f(t_n,x,Y^n,Z^n,\Gamma^n),
\end{aligned}\right.
\end{equation}
where $\bar Y^{n+j}$ and $\bar Z^{n+j}$
 are respectively the values of
$Y^{n+j}$ and $Z^{n+j}$ at the space points $\bar X^{n,j}$ that is defined by
\begin{equation}
\bar X^{n,j} = x + b(t_n,x)\Delta t_{n,j} + \sigma(t_n,x)\Delta W_{n,j},\quad j = 1,\ldots,k.\label{semi_sch_x2}
\end{equation}
Plug \eqref{semi_sch_x2} into the conditional expectation of  \eqref{semi_Eul}, since $\bar Y, \bar Z, \bar \Gm$
and $\bar A$ are all functions of $X^{n,j}$, we have
\begin{align*}
\Ec{t_n}{x}{\bar Y^{n+j}\Delta W_{n,j}^\tp} &=
 \Ec{t_n}{x}{Y^{n+j}(x + b(t_n,x)\Delta t_{n,j} + \sigma(t_n,x)\Delta W_{n,j})\Delta W_{n,j}^\tp}\\
 \vspace{0.35cm}
\Ec{t_n}{x}{\bar Z^{n+j}} &= \Ec{t_n}{x}{Z_{t_{n+j}}(x + b(t_n,x)\Delta t_{n,j} + \sigma(t_n,x)\Delta W_{n,j})}\\
\vspace{0.35cm}
\Ec{t_n}{x}{(\bar Z^{n+j})^\tp\Delta W_{n,j}^\tp}& =
\Ec{t_n}{x}{(Z^{n+j})^\tp(x + b(t_n,x)\Delta t_{n,j} + \sigma(t_n,x)\Delta W_{n,j})\Delta W_{n,j}^\tp}\\
\vspace{0.35cm}
\Ec{t_n}{x}{\bar Y^{n+j}} &= \Ec{t_n}{x}{ Y^{n+j}(x + b(t_n,x)\Delta t_{n,j} + \sigma(t_n,x)\Delta W_{n,j})}
\end{align*}
Note that quadrature methods should be applied when approximating the conditional
expectations above. Any efficient
quadrature rules such as the Monte-Carlo methods,
the quasi-Monte-Carlo methods, and the Gaussian quadrature methods can be used. Here, we will use the Gaussian quadratures based on the zeros of Hermite polynomials, for details, one can refer to \cite{ZFZ2014}.
In what follows, we shall denote by $\HE{\cdot}$ the numerical quadrature operator for
the conditional expectations.

\newcommand{\I}{\mathbb{I}}
Furthermore, when quadrature method is applied to approximate the conditional expectations, non-grid information might be needed. That's to say, for  $x\in \D_{h}^n$, points ($X^{n,j}$ defined by \eqref{semi_sch_x2} e.t.c.)
that are not in $\D_{h}^{n+j}$ are needed
when approximate the conditional expectation at time level $t_{n+j}$.
Thus, interpolation methods are also needed. We shall denote by $\I_\D$,
a local interpolation operator, such that $\I_\D g$ is the continuous function interpolated by the values of  $g$ at the grid points $\D$. We also denote by $\I_\D^n$ with the superscript $n$ indicating the
interpolation method at time level $t_n$. If the grid is a neighbor grid $U_{n,x}^{\D}$ of $x$, we denote by
$\I_{\D,x}^n$ instead of $\I_{U_{n,x}^{\D}}^n,$ for simplicity. Note that any interpolation methods can be used here,
however, care should be made if one wants to guarantee the stability and accuracy.

Now by introducing the operators
$\HE{\cdot}$ and $\I_{\D,x}^n$ (or more precisely
 $\I_{\D,\bar X_{t_{n+j}}^{t_n,x}}^{n+j}$)
one can rewrite \eqref{ref_all} in the following equivalent form
\begin{equation}\label{FAP_YZ}
\begin{aligned}
Z_{t_n} = & \sum\limits_{j=1}^k\alpha_{k,j}
\HE{\I_{\D,\bar X_{t_{n+j}}^{t_n,x}}^{n+j} Y_{t_{n+j}}\Delta W_{n,j}^\tp}
-R_{z,n}^k + R_{z,n}^{k,\mathbb{E}} + R_{z,n}^{k,\I},\\
\Gamma_{t_n} = & \sum\limits_{j=1}^k\alpha_{k,j}
\HE{\I_{\D,\bar X_{t_{n+j}}^{t_n,x}}^{n+j} Z_{t_{n+j}}^\tp\Delta W_{n,j}^\tp}
-R_{\Gamma,n}^k + R_{\Gamma,n}^{k,\mathbb{E}} + R_{\Gamma,n}^{k,\I},\\
A_{t_n} = & \sum\limits_{j=0}^k\alpha_{k,j}
\HE{\I_{\D,\bar X_{t_{n+j}}^{t_n,x}}^{n+j} Z_{t_{n+j}}}
-R_{A,n}^k + R_{A,n}^{k,\mathbb{E}} + R_{A,n}^{k,\I},\\
-\alpha_{k,0}Y_{t_n} =& \sum\limits_{j=1}^k\alpha_{k,j}
\HE{\I_{\D,\bar X_{t_{n+j}}^{t_n,x}}^{n+j} Y_{t_{n+j}}}
+f(t_n,x,Y_{t_n},Z_{t_n})\\
&+ R_{y,n}^k + R_{y,n}^{k,\mathbb{E}}+ R_{y,n}^{k,\I},
\end{aligned}
\end{equation}
where
$$\begin{array}{rl}
R_{z,n}^{k,\mathbb{E}} &= \sum_{j=1}^k\alpha_{k,j}(\mathbb E_{t_n}^{x}-\hat{\mathbb{E}}^{n,x})\left[{\bar Y_{t_{n+j}}\Delta W_{n,j}^\tp}\right],\\
R_{z,n}^{k,\I} &= \sum_{j=1}^k\alpha_{k,j}\Ec{t_n}{x,h}{(\bar Y_{t_{n+j}}-\I_{\D,\bar X_{t_{n+j}}^{t_n,x}}^{n+j} Y_{t_{n+j}})\Delta W_{n,j}^\tp},\\
R_{A,n}^{k,\mathbb{E}} &= \sum_{j=1}^k\alpha_{k,j}(\mathbb E_{t_n}^{x}-\hat{\mathbb{E}}^{n,x})\left[{\bar Z_{t_{n+j}}}\right],\\
R_{A,n}^{k,\I} &= \sum_{j=1}^k\alpha_{k,j}\Ec{t_n}{x,h}{(\bar Z_{t_{n+j}}-\I_{\D,\bar X_{t_{n+j}}^{t_n,x}}^{n+j} Z_{t_{n+j}})},\\
R_{\Gamma,n}^{k,\mathbb{E}} &= \sum_{j=1}^k\alpha_{k,j}(\mathbb E_{t_n}^{x}-\hat{\mathbb{E}}^{n,x})\left[{\bar Z_{t_{n+j}}^\tp\Delta W_{n,j}^\tp}\right],\\
R_{\Gamma,n}^{k,\I} &= \sum_{j=1}^k\alpha_{k,j}\Ec{t_n}{x,h}{(\bar Z_{t_{n+j}}^\tp-\I_{\D,\bar X_{t_{n+j}}^{t_n,x}}^{n+j} Z_{t_{n+j}}^\tp)\Delta W_{n,j}^\tp},\\
R_{y,n}^{k,\mathbb{E}} &= -\sum_{j=1}^k\alpha_{k,j}(\mathbb E_{t_n}^{x}-\hat{\mathbb{E}}^{n,x})\left[{\bar Y^{n+j}}\right],\\
R_{y,n}^{k,\I} &= -\sum_{j=1}^k\alpha_{k,j}\Ec{t_n}{x,h}{\bar Y_{t_{n+j}}-\I_{\D,\bar X_{t_{n+j}}^{t_n,x}}^{n+j} Y_{t_{n+j}}}.
\end{array}$$
The four terms $R_{y,n}^{k,\mathbb{E}}$, $R_{z,n}^{k,\mathbb{E}}$,
$R_{A,n}^{k,\mathbb{E}}$ and $R_{\Gamma,n}^{k,\mathbb{E}}$ are the error terms
resulted from approximating conditional expectations,
and the other four terms $R_{y,n}^{k,\I}$, $R_{z,n}^{k,\I}$,
$R_{A,n}^{k,\I}$ and $R_{\Gamma,n}^{k,\I}$ are the error terms
caused by numerical interpolations.

By removing those error terms $R_{y,n}^k$, $R_{y,n}^{k,\mathbb{E}}$, $ R_{y,n}^{k,\I}$,
$R_{z,n}^k$, $R_{z,n}^{k,\mathbb{E}}$, $R_{z,n}^{k,\I}$,
 $R_{A,n}^k$, $R_{A,n}^{k,\mathbb{E}}$, $ R_{A,n}^{k,\I}$,
$R_{\Gamma,n}^k$, $R_{\Gamma,n}^{k,\mathbb{E}}$ and $R_{\Gamma,n}^{k,\I}$
from \eqref{FAP_YZ},
we propose our fully discrete scheme for solving 2FBSDEs as follows:
\begin{scheme}\label{full_sch}
Assume random variables $Y^{N-i}$ and $Z^{N-i}$ defined on $\D_h^{N-i}$, $i = 0,1,\ldots,k-1$,
are known. For $n=N-k,\ldots,0$,
and for each $x\in \D_{h}^n$, solve $X^n$, $Y^n$, $Z^n$, $A^n$ and $\Gamma^n$ by
\begin{flalign}
X^{n,j} &= x + b(t_n,x)\Delta t_{n,j} + \sigma(t_n, x)\Delta W_{n,j},\quad j = 1,\ldots,k,\label{full_sch_x}
\end{flalign}\begin{flalign}
Z^n &= \sum_{j=1}^k\alpha_{k,j}
\HE{\I_{\D,\bar X^{n,j}}^{n+j} Y^{n+j}\Delta W_{n,j}^\tp},\label{full_sch_z}\\
\Gamma^n &= \sum_{j=1}^k\alpha_{k,j}
\HE{\I_{\D,\bar X^{n,j}}^{n+j} (Z^{n+j})^\tp\Delta W_{n,j}^\tp},\label{full_sch_Ga}\\
A^n &= \sum_{j=0}^k\alpha_{k,j}
\HE{\I_{\D,\bar X^{n,j}}^{n+j} Z^{n+j}}, \label{full_sch_A}\\
-\alpha_{k,0}Y^n &= \sum_{j=1}^k\alpha_{k,j}
\HE{\I_{\D,\bar X^{n,j}}^{n+j} Y^{n+j}}
+f(t_n,x,Y^n,Z^n,\Gamma^n).\label{full_sch_y}
\end{flalign}
\end{scheme}

In Scheme 3, for a fixed integer $k$, to get the values $Y^n$, $Z^n$, $A^n$ and $\Gamma^n$ on each time level $t_n$ at the grid $\D_h^n$,  we shall solve the values at $x \in \D_h^n$ one by one. Since they are independent with each other, the precess can be completely parallel. For a fixed $x\in \D_h^n$, in Scheme \ref{full_sch} , we firstly solve $X^{n,j}$ by the Euler scheme \eqref{full_sch_x} for $1\le j \le k$; and then solve $Z^n$, $\Gamma^n$ and $A^n$ by \eqref{full_sch_z}, \eqref{full_sch_Ga} and \eqref{full_sch_A} explicitly; finally, solve $Y^n$ by \eqref{full_sch_y} implicitly. Thus some iteration methods are required for solving $Y^n.$

If the function $f(t_n,x,y,z,\gamma)$
is Lipschitz continuous with respect to $y$, for small time partition step size $\Delta t_n$,
we can use the following
iteration procedure to solve $Y^n$
\begin{flalign}
\alpha_{k,0}Y^{n,l+1} &= -\sum_{j=1}^k\alpha_{k,j}
\HE{\I_{\D,\bar X^{n,j}}^{n+j} Y^{n+j}}
-f(t_n,x,Y^{n,l},Z^n,\Gamma^n),
\end{flalign}
with the iteration stopping condition $|Y^{n,l+1}-Y^{n,l}|\le \epsilon_0,$
where $\epsilon_0>0$ is a given tolerance. In case the function $f(t_n, x, y, z, \gamma)$ is differentiable with respect to $y$, to accelerate the convergence rate, the Newton iteration method can be applied , i.e.
\begin{equation*}
Y^{n,l+1}  = Y^{n,l} -\f{ \alpha_{k,0}Y^{n,l} +\sum\limits_{j=1}^k\alpha_{k,j}
\HE{\I_{\D,\bar X^{n,j}}^{n+j} Y^{n+j}}
+f(t_n,x,Y^{n,l},Z^n,\Gamma^n)}{\alpha_{k,0} + f_y(t_n,x,Y^{n,l} ,Z^n,\Gamma^n)}
\end{equation*}

Note that the local truncation errors of Scheme \ref{full_sch}
consist of twelve terms $R_{y,n}^k$, $R_{y,n}^{k,\mathbb{E}}$, $ R_{y,n}^{k,\I}$,
$R_{z,n}^k$, $R_{z,n}^{k,\mathbb{E}}$, $R_{z,n}^{k,\I}$,
 $R_{A,n}^k$, $R_{A,n}^{k,\mathbb{E}}$, $ R_{A,n}^{k,\I}$,
$R_{\Gamma,n}^k$, $R_{\Gamma,n}^{k,\mathbb{E}}$ and $R_{\Gamma,n}^{k,\I}$.
The four terms  $R_{y,n}^k$, $R_{z,n}^k$, $R_{A,n}^k$ and $R_{\Gamma,n}^k$
 defined respectively in \eqref{app}
come from the approximations of the derivative,
and the four terms $ R_{y,n}^{k,\I}$, $R_{z,n}^{k,\I}$,
$ R_{A,n}^{k,\I}$ and $R_{\Gamma,n}^{k,\I}$ defined in \eqref{FAP_YZ}
are the local interpolation errors.
For these eight terms, when data $b$, $\sigma$, $f$ and $g$
are smooth, the following estimates hold (provided that degree $r$ polynomials interpolation is used)
\begin{equation}\label{erro}
R_{y,n}^k, \,\, R_{z,n}^k, \,\,
R_{A,n}^k, \,\, R_{\Gamma,n}^k=O(\Delta t)^k), \,\,
R_{z,n}^{k,\I}, \,\, R_{y,n}^{k,\I},\,\,
R_{A,n}^{k,\I}, \,\, R_{\Gamma,n}^{k,\I}=O\left(h^{r+1}\right).
\end{equation}
The other four terms $R_{y,n}^{k,\mathbb{E}}$, $R_{z,n}^{k,\mathbb{E}}$,
$R_{A,n}^{k,\mathbb{E}}$ and $R_{\Gamma,n}^{k,\mathbb{E}}$
are the local truncation errors resulted from the approximations of
the conditional expectations. Noticed that
these conditional expectations
are functions of Gaussian random variables, thus, one can construct efficient Gauss-Hermite quadrature rule for their approximations.

\section{Extensions to coupled 2FBSDEs}

In this section we extend Scheme \ref{full_sch} to solve the fully coupled 2FBSDE. More precisely, we first update Scheme \ref{full_sch} into the following scheme

\begin{scheme}\label{full_nsch1}
Assume random variables $Y^{N-i}$  and  $Z^{N-i}$
defined on $\D_h^{N-i}$, $i = 0,1,\ldots,k-1$,
are known. For $n=N-k,\ldots,0$,
and for each $x\in \D_{h}^n$, solve $Y^n$, $Z^n$, $A^n$ and $\Gamma^n$ by
\begin{equation}
\left\{
\begin{aligned}
X^{n,j}  & =  x + b(t_n,x,Y^n,Z^n,\Gamma^n)\Delta t_{n,j}
+ \sigma(t_n,x,Y^n,Z^n,\Gamma^n)\Delta W_{n,j}, \\
 &\qquad  j = 1,2,\ldots,k, \\
Z^n  & =  \sum_{j=1}^k\alpha_{k,j}
\HE{\I_{\D,\bar X^{n,j}}^{n+j} Y^{n+j}\Delta W_{n,j}^\tp},\\
\Gamma^n  & =  \sum_{j=1}^k\alpha_{k,j}
\HE{\I_{\D,\bar X^{n,j}}^{n+j} (Z^{n+j})^\tp\Delta W_{n,j}^\tp},\\
A^n  & =  \sum_{j=0}^k\alpha_{k,j}
\HE{\I_{\D,\bar X^{n,j}}^{n+j} Z^{n+j}},\\
-\alpha_{k,0}Y^n  & =  \sum_{j=1}^k\alpha_{k,j}
\HE{\I_{\D,\bar X^{n,j}}^{n+j} Y^{n+j}}
+f(t_n,x,Y^n,Z^n,\Gamma^n).
\end{aligned}
\right.
\end{equation}
\end{scheme}

The main difference in Scheme \ref{full_nsch1} is that $X^{n,j}$, $Y^n$,$Z^n$ and $\Gm^n$ are
all coupled together. Thus, one needs to solve nonlinear
equations with some iterative procedure. In our numerical experiments,
we use the following iterative scheme.

\begin{scheme}\label{full_nsch}
Assume random variables $Y^{N-i}$, $Z^{N-i}$ and $\Gamma^{N-i}$
defined on $\D_h^{N-i}$, $i = 0,1,\ldots,k-1$,
are known. $l$ is  the iterative time,  $Y^{n,l}$, $Z^{n,l}$, $A^{n,l}$ and $\Gamma^{n,l}$ are the corresponding values at the $l$-th iterative.
For $n=N-k,\ldots,0$, and  each $x\in \D_{h}^n$,  solve $Y^n$, $Z^n$, $A^n$ and $\Gamma^n$ by
\begin{enumerate}
\item Initialize $Y^{n,0}$, $Z^{n,0}$ and  $\Gamma^{n,0}$
by $Y^{n+1}$, $Z^{n+1}$ and $\Gamma^{n+1}$, separately,
i.e.  $Y^{n,0}= Y^{n+1}(x)$, $Z^{n,0}=Z^{n+1}(x)$  and $\Gamma^{n,0}=\Gm^{n+1}(x)$;
\item  Calculate $Y^{n,l+1}$,  $Z^{n,l+1}$,  $A^{n,l+1}$ and $\Gamma^{n,l+1}$  for $l = 0,1,\dots$ by
\begin{equation}
\left\{
\begin{aligned}
X^{n,j} = x& + b(t_n,x,Y^{n,l},Z^{n,l},\Gm^{n,l})\Delta t_{n,j} + \sigma(t_n,x,Y^{n,l},Z^{n,l},\Gm^{n,l})\Delta W_{n,j},\\
&\qquad j = 1,2,\ldots,k, \\
Z^{n,l+1} &= \sum_{j=1}^k\alpha_{k,j}\HE{\I_{\D,\bar X^{n,j}}^{n+j} Y^{n+j}\Delta W_{n,j}^\tp},\\
\Gamma^{n,l+1} &= \sum_{j=1}^k\alpha_{k,j}\HE{\I_{\D,\bar X^{n,j}}^{n+j} (Z^{n+j})^\tp\Delta W_{n,j}^\tp},\\
 A^{n,l+1} &= \sum_{j=0}^k\alpha_{k,j}\HE{\I_{\D,\bar X^{n,j}}^{n+j} Z^{n+j}},\nonumber\\
-\alpha_{k,0}Y^{n,l+1} &= \sum_{j=1}^k\alpha_{k,j}\HE{\I_{\D,\bar X^{n,j}}^{n+j} Y^{n+j}}
+f(t_n,x,Y^{n,l+1},Z^{n,l+1},\Gm^{n,l+1}).
\end{aligned}\right.
\end{equation}
Repeat the above calculation, until
$$\max\left\{|Y^{n,l+1}-Y^{n,l}|,|Z^{n,l+1}-Z^{n,l}|,
|A^{n,l+1}-A^{n,l}|,|\Gamma^{n,l+1}-\Gamma^{n,l}|\right\}< \epsilon_0;$$
\item The values $(Y^n, Z^n, \Gm^n,A^n)$ are set to be $(Y^{n,l+1}, Z^{n,l+1},
 A^{n,l+1}, \Gamma^{n,l+1})$.
\end{enumerate}
The computation is complete.
\end{scheme}

It is worth to remark that our Schemes \ref{full_nsch1} and \ref{full_nsch}
are heuristic generalizations of
the decoupled ones in the last section.
In case the drift coefficient $b$ and the diffusion
coefficient $\sigma$ do not depend on $Y$, $Z$ and $\Gamma$,
Scheme \ref{full_nsch} coincides with Scheme \ref{full_sch}.
\begin{remark}
We have finished the construction of our multistep numerical schemes, i.e., Scheme 3 for decoupled 2FBSDEs and Scheme 5 for coupled 2FBSDEs. In both cases, the Euler method is used to discrete the forward SDEs, thus, the total computational complexity can be significantly reduced. Furthermore, we shall show by numerical tests that the solution $(Y_t,Z_t,\Gamma_t,A_t)$ can still be of high order accuracy, in the next section.
\end{remark}
\section{Numerical experiments}\label{NuEx}

In this section, we shall provide with several constructive numerical tests to show the efficiency and accuracy of our multistep schemes proposed in the previous sections. Applications of our numerical methods for stochastic optimal control problems will also be presented.

In all our numerical tests, we shall consider the uniform partition for the time-space domain $[0,T]\times \R^m$, for simplicity. That is, for a given $N>0, h>0$, we set
 \begin{align*}
 \mathcal{T} &:= \{t_n\ | \ t_n = n \Delta t, \quad n=0,1,\dots,N, \Delta t = \f {T}N\}, \\
 \D_h & := \{x_\bfj\ |\ x_\bfj = x_0 + \bfj \cdot h,\quad   \bfj :=(j_1,j_2,\dots,j_m)^\tp, \text{for each}\ j_i \in \mathbf  Z\},
 \end{align*}
with $x_0\in \R^m$ the initial state. Specially, for the one dimensional case, i.e.\  $m=1$, this yields $\D_h = \{x_j \ |\  x_j = x_0 + jh, j=0,\pm1,\pm2,\dots\}$. Our scheme is applied to calculate
the values of  $Y^n$, $Z^n$, $\Gm^n$ and $A^n$,
which are the numerical values of $Y$, $Z$, $\Gm$ and $A$
at every time-space point $(t_n, x_\bfj) \in \mathcal{T}\times \D_h$.

In our tests, we shall adopt the local Lagrange interpolation method $\I_\D$, so that the interpolation error estimates in \eqref{erro} holds, and the Hermite-Gaussian quadrature rule will be used to approximate the associated conditional expectations. Since we aim to checking the convergence rates of our multistep schemes, all errors
resulted from these approximations are controlled by choosing
appropriate parameters (e.g. interpolation order $r$ and the numerical of Gaussian-Hermite points). In particular, we shall use 10 Gauss-Hermite quadrature points in each dimension, and for the Lagrange interpolation, more than 6 points will be used. To balance the errors result from the time discrete truncation  and the space truncation, we choose $h=(\Delta t)^{\frac{k+1}{r+1}},$ where $r$ is the degree of the Langrange interpolation polynomial. Furthermore, for the $k$-step schemes, the information for $\{Y^{N-j},Z^{N-j}\}_{j=1}^k$ are needed, this can be obtained by using other numerical methods with small time steps, to maintain the high order convergence rates. Here, we just artificially assume that the values are known for simplicity.   

The numerical results are obtained with \texttt{FORTRAN 95}
on a workstation with one \texttt{Intel Xeon E5-2620 v2 CPU (12 cores, 2.10 GHz )}.
To accelerate the performance, \texttt{OpenMP} techniques are used.
To guarantee the computing precision, we use the \textit{long double} type (\texttt{real(16)}) digital for the float variables when programming. The long double variable has 34 significant digits which supplies enough
precision for our computation. However, the cost of time is increasing dramatically compared with variables all
defined as double (\texttt{real(8)}), but even so, the time elapsed of our programs is still competitive.
In what follows, we will denote by CR the convergence rates and  $T_r$ the running time, respectively. For all our numerical examples, the terminal time $T$ is set to be $1.0.$

\subsection{Decoupled and coupled 2FBSDEs}
We first test Scheme 3 for solving decoupled 2FBSDEs. The first example considered is
\begin{equation}
\left\{\begin{aligned}
\d X_t &=  \sin(t+X_t) \d t + c \cos(t+ X_t) \d W_t,  \\
-\d Y_t &= \big(-\cos(t+X_t)\f1cZ_t-\cos(t+X_t)(Y_t^2+Y_t) -\f14\Gamma_t \big)\d t - Z_t \d W_t,\\
\d Z_t &= A_t \d t + \Gamma_t \d W_t,
\end{aligned}\right.
\end{equation}
with the initial condition $X_{t_0} = x$ and the terminal condition
$Y_T = \sin(T+X_T).$ It can be shown that the exact solutions are
\begin{align*}
&Y_t= \sin(t+X_t), \quad Z_t=c \cos^2(t+X_t), \quad \Gamma_t=  -2c^2\sin(t+X_t)\cos^2(t+X_t),\\
&A_t= -c\sin(2t+2X_t)(1+\sin(t+X_t)) -c^3\cos(2t+2X_t)\cos^2(t+X_t).
\end{align*}
In the numerical test, we set $x = 0.5$ and $t_{0}=0$, and solve the 2FBSDEs by scheme 3 with different step parameter $k.$ The numerical errors $\abs{Y^0-Y_0}, \abs{Z^0-Z_0},\abs{\Gm^0-\Gm_0}$ and $\abs{A^0-A_0}$
and the corresponding convergence rates are listed in Table \ref{tab:1}.
It is shown in Table \ref{tab:1} that: the multistep numerical scheme works very well, and the numerical error goes to machine accuracy. Moreover, the method admits a $k$-order convergence rate, and it remains stable for $1\le k \le 6,$ which is coincide with the classic numerical ODEs theory and our previous results \cite{ZFZ2014}.

Furthermore, given a fixed accuracy tolerance, it would be more efficient if the multistep scheme with a large $k$ is used.   This can be more easily seen from the following table.

\vspace{1em}
\begin{tabular}{*{6}{c|}c}
\hline
STEP & N & $T_{rn}$  &  $\abs{Y^0-Y_0}$   &  $\abs{Z^0-Z_0}$  &  $\abs{\Gm^0-\Gm_0}$  & $\abs{A^0-A_0}$ \\
\hline
$k = 1$  & 8192 & 135.0s & 5.172E-05 & 2.739E-05 & 7.573E-07 &  3.585E-07 \\
\hline
$k=2 $  & 2048 & 10.82s & 1.232E-09 &  4.704E-08 &  2.301E-08 &  1.920E-07 \\
\hline
$k=3$  &  512    &  5.58s   & 3.423E-09 &  3.097E-09 &  7.797E-10 &  4.677E-09 \\
\hline
$k=4$  &  128    & 2.85s     &  1.159E-09 &  1.733E-09 &  2.113E-09 &  6.042E-09 \\
\hline
\end{tabular}
\vspace{1em}
\begin{table}[!htbp]
\caption{Numerical results of example 1 with $c=0.1$}\label{tab:1}
\begin{minipage}{\textwidth}\footnotesize%
\begin{center}%
\begin{tabular}{*{7}{c|}c}\hline
STEP& &N  &  $\abs{Y^0-Y_0}$   &  $\abs{Z^0-Z_0}$  &  $\abs{\Gamma^0-\Gamma_0}$  & $\abs{A^0-A_0}$  &  $T_{\mathrm rn}$ \\
\hline
\rb{-3em}{K=1} & \rb{-2.5em}{DT}
&   32 &  1.253E-02 &  6.800E-03 &  3.156E-04 &  1.703E-03 &
\\ \cline{3-7} &
&   64 &  6.411E-03 &  3.466E-03 &  1.854E-04 &  1.059E-03 &
\\ \cline{3-7} &
&  128 &  3.248E-03 &  1.752E-03 &  8.667E-05 &  4.409E-04 &
\\ \cline{3-7} &
&  256 &  1.637E-03 &  8.808E-04 &  3.570E-05 &  1.381E-04 &
\\ \cline{3-7} &
&  512 &  8.227E-04 &  4.405E-04 &  1.411E-05 &  3.472E-05 &
\\ \cline{2-7} & CR &&  0.98&  0.99&  1.13&  1.42& \rb{3em}{  1.87s}
\footnote{\scriptsize -k 1 -ng 10 -ni 5 -e  F -sh  1.0 -mg -20.00 20.00 -$t_0$0.0 -$x_0$0.5 -T1.0}
\\ \hline
\rb{-3em}{K=2} & \rb{-2.5em}{DT}
&   32 &  6.123E-06 &  1.897E-04 &  1.086E-04 &  9.312E-04 &
\\ \cline{3-7} &
&   64 &  1.334E-06 &  4.761E-05 &  2.614E-05 &  2.210E-04 &
\\ \cline{3-7} &
&  128 &  3.332E-07 &  1.192E-05 &  6.207E-06 &  5.196E-05 &
\\ \cline{3-7} &
&  256 &  8.494E-08 &  2.995E-06 &  1.500E-06 &  1.253E-05 &
\\ \cline{3-7} &
&  512 &  2.133E-08 &  7.514E-07 &  3.707E-07 &  3.098E-06 &
\\ \cline{2-7} & CR &&  2.03&  2.00&  2.05&  2.06& \rb{3em}{  4.96s}
\footnote{\scriptsize -k 2 -ng 10 -ni 5 -e  F -sh  1.0 -mg -20.00 20.00 -$t_0$0.0 -$x_0$0.5 -T1.0}
\\ \hline
\rb{-3em}{K=3} & \rb{-2.5em}{DT}
&   32 &  1.149E-05 &  1.409E-05 &  2.846E-06 &  1.988E-05 &
\\ \cline{3-7} &
&   64 &  1.615E-06 &  1.668E-06 &  3.734E-07 &  2.382E-06 &
\\ \cline{3-7} &
&  128 &  2.117E-07 &  2.026E-07 &  4.848E-08 &  2.972E-07 &
\\ \cline{3-7} &
&  256 &  2.704E-08 &  2.496E-08 &  6.202E-09 &  3.727E-08 &
\\ \cline{3-7} &
&  512 &  3.415E-09 &  3.096E-09 &  7.818E-10 &  4.650E-09 &
\\ \cline{2-7} & CR &&  2.93&  3.04&  2.96&  3.01& \rb{3em}{ 13.96s}
\footnote{\scriptsize -k 3 -ng 10 -ni 5 -e  F -sh  1.0 -mg -20.00 20.00 -$t_0$0.0 -$x_0$0.5 -T1.0}
\\ \hline
\rb{-3em}{K=4} & \rb{-2.5em}{DT}
&   32 &  2.853E-07 &  4.169E-07 &  6.203E-07 &  2.171E-06 &
\\ \cline{3-7} &
&   64 &  1.843E-08 &  2.700E-08 &  3.588E-08 &  1.122E-07 &
\\ \cline{3-7} &
&  128 &  1.161E-09 &  1.736E-09 &  2.097E-09 &  5.886E-09 &
\\ \cline{3-7} &
&  256 &  7.271E-11 &  1.106E-10 &  1.275E-10 &  3.463E-10 &
\\ \cline{3-7} &
&  512 &  4.545E-12 &  6.973E-12 &  7.852E-12 &  2.097E-11 &
\\ \cline{2-7} & CR &&  3.99&  3.97&  4.07&  4.17& \rb{3em}{ 17.74s}
\footnote{\scriptsize -k 4 -ng 10 -ni 8 -e  F -sh  1.0 -mg -20.00 20.00 -$t_0$0.0 -$x_0$0.5 -T1.0}
\\ \hline
\rb{-3em}{K=5} & \rb{-2.5em}{DT}
&   32 &  1.124E-08 &  3.151E-08 &  2.563E-08 &  1.660E-07 &
\\ \cline{3-7} &
&   64 &  5.297E-10 &  9.433E-10 &  8.493E-10 &  4.027E-09 &
\\ \cline{3-7} &
&  128 &  1.877E-11 &  2.770E-11 &  2.632E-11 &  9.995E-11 &
\\ \cline{3-7} &
&  256 &  6.180E-13 &  8.174E-13 &  8.310E-13 &  2.951E-12 &
\\ \cline{3-7} &
&  512 &  1.971E-14 &  2.480E-14 &  2.693E-14 &  9.838E-14 &
\\ \cline{2-7} & CR &&  4.80&  5.07&  4.97&  5.18& \rb{3em}{ 34.41s}
\footnote{\scriptsize -k 5 -ng 10 -ni 14 -e  F -sh  1.0 -mg -20.00 20.00 -$t_0$0.0 -$x_0$0.5 -T1.0}
\\ \hline
\rb{-3em}{K=6} & \rb{-2.5em}{DT}
&   32 &  4.721E-10 &  9.254E-10 &  5.225E-09 &  1.144E-08 &
\\ \cline{3-7} &
&   64 &  7.847E-12 &  1.542E-11 &  7.602E-11 &  1.545E-10 &
\\ \cline{3-7} &
&  128 &  1.144E-13 &  2.506E-13 &  1.140E-12 &  2.144E-12 &
\\ \cline{3-7} &
&  256 &  1.741E-15 &  4.074E-15 &  1.710E-14 &  2.863E-14 &
\\ \cline{3-7} &
&  512 &  2.674E-17 &  6.557E-17 &  2.622E-16 &  4.232E-16 &
\\ \cline{2-7} & CR &&  6.03&  5.94&  6.06&  6.18& \rb{3em}{ 50.75s}
\footnote{\scriptsize -k 6 -ng 10 -ni 16 -e  F -sh  1.0 -mg -20.00 20.00 -$t_0$0.0 -$x_0$0.5 -T1.0}
\\ \hline
\rb{-3em}{K=7} & \rb{-2.5em}{DT}
&   32 &  1.422E-11 &  1.805E-11 &  4.311E-10 &  7.197E-09 &
\\ \cline{3-7} &
&   64 &  3.876E-13 &  8.371E-13 &  1.417E-12 &  3.941E-12 &
\\ \cline{3-7} &
&  128 &  2.365E-15 &  1.959E-15 &  2.089E-14 &  3.823E-13 &
\\ \cline{3-7} &
&  256 &  2.762E-17 &  1.366E-17 &  1.265E-15 &  1.528E-14 &
\\ \cline{3-7} &
&  512 &  4.993E-18 &  1.958E-16 &  9.029E-17 &  1.894E-16 &
\\ \cline{2-7} & CR &&  5.67&  4.89&  5.45&  5.84& \rb{3em}{ 83.65s}
\footnote{\scriptsize -k 7 -ng 10 -ni 20 -e  F -sh  1.0 -mg -20.00 20.00 -$t_0$0.0 -$x_0$0.5 -T1.0}
\\ \hline
\rb{-3em}{K=8} & \rb{-2.5em}{DT}
&   32 &  7.015E-12 &  1.501E-10 &  2.182E-09 &  2.147E-08 &
\\ \cline{3-7} &
&   64 &  4.860E-12 &  6.270E-11 &  3.480E-10 &  7.047E-09 &
\\ \cline{3-7} &
&  128 &  1.117E-10 &  6.105E-09 &  8.660E-08 &  1.014E-06 &
\\ \cline{3-7} &
&  256 &  7.367E-03 &  3.798E+00 &  4.029E+02 &  3.354E+04 &
\\ \cline{3-7} &
&  512 &        NaN &        NaN &        NaN &        NaN &
\\ \cline{2-7} & CR &&   NaN&   NaN&   NaN&   NaN& \rb{3em}{108.42s}
\footnote{\scriptsize -k 8 -ng 10 -ni 20 -e  F -sh  1.0 -mg -20.00 20.00 -$t_0$0.0 -$x_0$0.5 -T1.0}
\\ \hline
\end{tabular}\end{center}\end{minipage}
\end{table}

We now consider an example with the geometry Brownian motion in the forward SDE, the main feature here is that the drift and diffusion terms are unbounded. The example yields
\begin{eqsys}
\d X_t &= r X_t \d t + c X_t \d W_t,\\
-\d Y_t&= -e^{-\f{X_t^2}M} +\f2MrX_t^2Y_t-\f{\Gm_t}2+\f c2Z_t \d t -Z_t \d W_t,\\
\d Z_t&= A_t \d t + \Gm_t \d W_t
\end{eqsys}
with terminal condition $Y_T = Te^{-\f{X_T^2}M}.$ It can be shown that the exact solution is
\begin{align*}
&Y_t = te^{-\f{X_t^2}M}, \quad Z_t= -\f{2c}MtX_t^2e^{-\f{X_t^2}M}, \quad  \Gm_t= \f{4c^2}{M^2}tX_t^2e^{-\f{X_t^2}M}(X_t^2 -M),\\
&A_t = -\f{2c}Me^{-\f{X_t^2}M}[(1+2rt+c^2t)X_t^2 -\f{2r+c^2}MtX_t^4 + \f{2c^2t}MX_t^6].
\end{align*}
In this example, we set $x_0=1.5,$ and the parameters are chosen as $r = 0.2, c=0.01, M=4.$ The corresponding numerical results are listed the in Table \ref{eg:geo}. Again, our multistep schemes behaves very well, and high order convergence rates are obtained.
\begin{table}[!hbp]
\caption{Numerical results of example 2.}\label{eg:geo}
\begin{minipage}{\textwidth}\footnotesize%
\begin{center}%
\begin{tabular}{*{7}{c|}c}\hline
STEP& &N  &  $\abs{Y^0-Y_0}$   &  $\abs{Z^0-Z_0}$  &  $\abs{\Gamma^0-\Gamma_0}$  & $\abs{A^0-A_0}$  &  $T_{\mathrm rn}$ \\
\hline
\rb{-3em}{K=1} & \rb{-2.5em}{DT}
&   32 &  4.808E-03 &  2.331E-04 &  2.431E-06 &  6.699E-07 &
\\ \cline{3-7} &
&   64 &  2.402E-03 &  1.169E-04 &  1.232E-06 &  3.335E-07 &
\\ \cline{3-7} &
&  128 &  1.200E-03 &  5.854E-05 &  6.207E-07 &  1.665E-07 &
\\ \cline{3-7} &
&  256 &  5.999E-04 &  2.929E-05 &  3.115E-07 &  8.308E-08 &
\\ \cline{3-7} &
&  512 &  2.999E-04 &  1.465E-05 &  1.560E-07 &  4.143E-08 &
\\ \cline{2-7} & CR &&  1.00&  1.00&  0.99&  1.00& \rb{3em}{  2.58s}
\footnote{\scriptsize -k 1 -ng 10 -ni 10 -e  F -sh  1.0 -mg -20.00 20.00 -$t_0$0.0 -$x_0$1.5 -T1.0}
\\ \hline
\rb{-3em}{K=2} & \rb{-2.5em}{DT}
&   32 &  1.385E-05 &  1.294E-06 &  1.459E-07 &  1.457E-06 &
\\ \cline{3-7} &
&   64 &  3.582E-06 &  3.252E-07 &  3.653E-08 &  3.603E-07 &
\\ \cline{3-7} &
&  128 &  9.098E-07 &  8.147E-08 &  9.140E-09 &  8.954E-08 &
\\ \cline{3-7} &
&  256 &  2.291E-07 &  2.039E-08 &  2.286E-09 &  2.231E-08 &
\\ \cline{3-7} &
&  512 &  5.747E-08 &  5.098E-09 &  5.717E-10 &  5.569E-09 &
\\ \cline{2-7} & CR &&  1.98&  2.00&  2.00&  2.01& \rb{3em}{  7.54s}
\footnote{\scriptsize -k 2 -ng 10 -ni 10 -e  F -sh  1.0 -mg -20.00 20.00 -$t_0$0.0 -$x_0$1.5 -T1.0}
\\ \hline
\rb{-3em}{K=3} & \rb{-2.5em}{DT}
&   32 &  3.490E-07 &  5.124E-08 &  1.213E-09 &  8.537E-09 &
\\ \cline{3-7} &
&   64 &  4.638E-08 &  6.484E-09 &  1.552E-10 &  1.130E-09 &
\\ \cline{3-7} &
&  128 &  5.983E-09 &  8.157E-10 &  1.960E-11 &  1.454E-10 &
\\ \cline{3-7} &
&  256 &  7.604E-10 &  1.023E-10 &  2.461E-12 &  1.846E-11 &
\\ \cline{3-7} &
&  512 &  9.590E-11 &  1.281E-11 &  3.082E-13 &  2.326E-12 &
\\ \cline{2-7} & CR &&  2.96&  2.99&  2.99&  2.96& \rb{3em}{ 16.73s}
\footnote{\scriptsize -k 3 -ng 10 -ni 10 -e  F -sh  1.0 -mg -20.00 20.00 -$t_0$0.0 -$x_0$1.5 -T1.0}
\\ \hline
\rb{-3em}{K=4} & \rb{-2.5em}{DT}
&   32 &  4.889E-09 &  8.085E-10 &  7.806E-11 &  6.895E-10 &
\\ \cline{3-7} &
&   64 &  3.369E-10 &  5.171E-11 &  4.934E-12 &  4.232E-11 &
\\ \cline{3-7} &
&  128 &  2.212E-11 &  3.276E-12 &  3.100E-13 &  2.601E-12 &
\\ \cline{3-7} &
&  256 &  1.417E-12 &  2.063E-13 &  1.942E-14 &  1.610E-13 &
\\ \cline{3-7} &
&  512 &  8.969E-14 &  1.295E-14 &  1.215E-15 &  1.003E-14 &
\\ \cline{2-7} & CR &&  3.94&  3.98&  3.99&  4.02& \rb{3em}{ 34.39s}
\footnote{\scriptsize -k 4 -ng 10 -ni 13 -e  F -sh  1.0 -mg -20.00 20.00 -$t_0$0.0 -$x_0$1.5 -T1.0}
\\ \hline
\rb{-3em}{K=5} & \rb{-2.5em}{DT}
&   32 &  4.916E-11 &  7.997E-12 &  6.357E-14 &  2.573E-12 &
\\ \cline{3-7} &
&   64 &  1.548E-12 &  2.514E-13 &  1.625E-15 &  7.997E-14 &
\\ \cline{3-7} &
&  128 &  4.545E-14 &  7.896E-15 &  3.486E-17 &  2.277E-15 &
\\ \cline{3-7} &
&  256 &  1.316E-15 &  2.481E-16 &  6.573E-19 &  6.360E-17 &
\\ \cline{3-7} &
&  512 &  3.843E-17 &  7.792E-18 &  1.067E-20 &  1.790E-18 &
\\ \cline{2-7} & CR &&  5.08&  4.99&  5.63&  5.12& \rb{3em}{ 61.72s}
\footnote{\scriptsize -k 5 -ng 10 -ni 16 -e  F -sh  1.0 -mg -20.00 20.00 -$t_0$0.0 -$x_0$1.5 -T1.0}
\\ \hline
\rb{-3em}{K=6} & \rb{-2.5em}{DT}
&   32 &  2.246E-12 &  3.449E-13 &  3.677E-14 &  4.326E-13 &
\\ \cline{3-7} &
&   64 &  3.805E-14 &  5.534E-15 &  5.708E-16 &  6.289E-15 &
\\ \cline{3-7} &
&  128 &  5.870E-16 &  8.807E-17 &  9.179E-18 &  9.283E-17 &
\\ \cline{3-7} &
&  256 &  8.808E-18 &  1.394E-18 &  1.463E-19 &  1.391E-18 &
\\ \cline{3-7} &
&  512 &  1.319E-19 &  2.197E-20 &  2.318E-21 &  2.136E-20 &
\\ \cline{2-7} & CR &&  6.01&  5.98&  5.98&  6.07& \rb{3em}{ 93.42s}
\footnote{\scriptsize -k 6 -ng 10 -ni 16 -e  F -sh  1.0 -mg -20.00 20.00 -$t_0$0.0 -$x_0$1.5 -T1.0}
\\ \hline
\rb{-3em}{K=7} & \rb{-2.5em}{DT}
&   32 &  3.908E-13 &  1.027E-13 &  2.830E-15 &  1.527E-13 &
\\ \cline{3-7} &
&   64 &  4.974E-16 &  1.217E-15 &  3.787E-17 &  1.929E-15 &
\\ \cline{3-7} &
&  128 &  3.750E-17 &  1.204E-17 &  3.601E-18 &  1.489E-16 &
\\ \cline{3-7} &
&  256 &  3.893E-19 &  2.617E-17 &  2.460E-15 &  3.265E-13 &
\\ \cline{3-7} &
&  512 &  5.287E-17 &  8.003E-15 &  6.419E-12 &  1.677E-09 &
\\ \cline{2-7} & CR &&  3.60&  1.29& -2.83& -3.42& \rb{3em}{196.54s}
\footnote{\scriptsize -k 7 -ng 10 -ni 30 -e  F -sh  1.0 -mg -20.00 20.00 -$t_0$0.0 -$x_0$1.5 -T1.0}
\\ \hline
\end{tabular}\end{center}\end{minipage}
\end{table}

We now test Scheme 5 for solving coupled 2FBSDEs. Here, iterative process is needed. To show the high accuracy of the proposed schemes, we set the tolerance as $\epsilon = 10^{-25}.$ We first consider the following coupled 2FBSDEs:
\begin{equation}
\left\{\begin{aligned}
\d X_t &=  (\sin(t+X_t) + Z_t/c +\sin(t+X_t)Y_t-1) \d t \\
&+ (c\cos(t+X_t)-c + c\cos^2(t+X_t) + c\sin(t+X_t)Y_t) \d W_t,\\
-\d Y_t &= -\cos(t+X_t)\f1cZ_t-\cos(t+X_t)(Y_t^2+Y_t) -\f14\Gm_t \d t - Z_t \d W_t,\\
\d Z_t &= A_t \d t + \Gm_t \d W_t,
\end{aligned}\right.
\end{equation}
with conditions $Y_T= \sin(T+X_T)$ and  $X_{t_0} = x.$ The exact solution is
\begin{align*}
&Y_t= \sin(t+X_t), \quad Z_t=c \cos^2(t+X_t), \quad \Gm_t=  -2c^2\sin(t+X_t)\cos^2(t+X_t),\\
&A_t= -c\sin(2t+2X_t)(1+\sin(t+X_t)) -c^3\cos(2t+2X_t)\cos^2(t+X_t).
\end{align*}
The numerical results are shown in Table \ref{dp:1}, where one can concludes that the numerical methods work well and  high-order convergence rates are obtained. However, due to the computational complexity, numerical results with only multistep up to $k=3$ are presented.

\begin{table}[!htbp]
\caption{Numerical results of example 3.}\label{dp:1}
\begin{minipage}{\textwidth}\footnotesize%
\begin{center}%
\begin{tabular}{*{7}{c|}c}\hline
STEP& &N  &  $\abs{Y^0-Y_0}$   &  $\abs{Z^0-Z_0}$  &  $\abs{\Gamma^0-\Gamma_0}$  & $\abs{A^0-A_0}$  &  $T_{\mathrm rn}$ \\
\hline

\rb{-3em}{K=1} & \rb{-2.5em}{DT}
&   32 &  4.126E-02 &  9.266E-03 &  3.980E-03 &  5.001E-02 &
\\ \cline{3-7} &
&   64 &  2.085E-02 &  4.782E-03 &  1.996E-03 &  2.532E-02 &
\\ \cline{3-7} &
&  128 &  1.048E-02 &  2.430E-03 &  9.968E-04 &  1.273E-02 &
\\ \cline{3-7} &
&  256 &  5.255E-03 &  1.225E-03 &  4.977E-04 &  6.377E-03 &
\\ \cline{3-7} &
&  512 &  2.631E-03 &  6.149E-04 &  2.486E-04 &  3.192E-03 &
\\ \cline{2-7} & CR &&  0.99&  0.98&  1.00&  0.99& \rb{3em}{ 35.42s}
\footnote{\scriptsize -k 1 -ng 10 -ni 5 -e  F -sh  1.0 -mg -20.00 20.00 -$t_0$0.0 -$x_0$1.0 -T1.0}
\\ \hline
\rb{-3em}{K=2} & \rb{-2.5em}{DT}
&   32 &  5.077E-04 &  4.045E-05 &  2.561E-04 &  1.640E-03 &
\\ \cline{3-7} &
&   64 &  1.250E-04 &  7.443E-06 &  6.648E-05 &  4.266E-04 &
\\ \cline{3-7} &
&  128 &  3.099E-05 &  1.541E-06 &  1.690E-05 &  1.086E-04 &
\\ \cline{3-7} &
&  256 &  7.713E-06 &  3.462E-07 &  4.258E-06 &  2.740E-05 &
\\ \cline{3-7} &
&  512 &  1.924E-06 &  8.166E-08 &  1.069E-06 &  6.876E-06 &
\\ \cline{2-7} & CR &&  2.01&  2.23&  1.98&  1.98& \rb{3em}{232.23s}
\footnote{\scriptsize -k 2 -ng 10 -ni 5 -e  F -sh  1.0 -mg -20.00 20.00 -$t_0$0.0 -$x_0$1.0 -T1.0}
\\ \hline
\rb{-3em}{K=3} & \rb{-2.5em}{DT}
&   32 &  6.995E-05 &  2.551E-05 &  2.388E-05 &  1.568E-04 &
\\ \cline{3-7} &
&   64 &  8.908E-06 &  3.156E-06 &  2.693E-06 &  1.840E-05 &
\\ \cline{3-7} &
&  128 &  1.123E-06 &  3.923E-07 &  3.173E-07 &  2.225E-06 &
\\ \cline{3-7} &
&  256 &  1.409E-07 &  4.889E-08 &  3.843E-08 &  2.733E-07 &
\\ \cline{3-7} &
&  512 &  1.764E-08 &  6.103E-09 &  4.726E-09 &  3.383E-08 &
\\ \cline{2-7} & CR &&  2.99&  3.01&  3.07&  3.04& \rb{3em}{987s}
\footnote{\scriptsize -k 3 -ng 10 -ni 5 -e  F -sh  1.0 -mg -20.00 20.00 -$t_0$0.0 -$x_0$1.0 -T1.0}
\\ \hline

\end{tabular}\end{center}\end{minipage}
\end{table}

Our next couple 2FBSDEs example yields
\begin{equation*}
\left\{\begin{aligned}
X_t &= x + \int_0^t \frac{1}{(1+\exp{(s+X_s)})(1+Y_s)} \di s + \int_0^t Y_s \di W_s,\\
Y_t &= \frac{\exp{(T+X_T)}}{1+\exp{(T+X_T)}}\\
& + \int_t^T \left[\frac{2Y_s}{1+2\exp{(s+X_s)}}
+ \frac12\left(\Gamma_s-\frac{Y_sZ_s}{1+\exp{(s+X_s)}}\right)\right] \di s - \int_t^T Z_s \di W_s,\\
Z_t &= Z_T - \int_t^T A_s \di s - \int_t^T \Gamma_s \di W_s.
\end{aligned}\right.
\end{equation*}
The associated exact solution is
\begin{align*}
&Y_t = \frac{\exp{(t+X_t)}}{1+\exp{(t+X_t)}}, \quad Z_t = \frac{\exp{(t+X_t)}^2}{(1+\exp{(t+X_t)})^3},\\
&\Gamma_t = \frac{\exp{(t+X_t)}^3(2-\exp{(t+X_t)})}{(1+\exp{(t+X_t)})^5},\\
&A_t = \frac{2\exp{(t+X_t)}^2(2-\exp{(t+X_t)})}{(1+\exp{(t+X_t)})^3(1+2\exp{(t+X_t)})}\\
&\qquad  + \frac{\exp{(t+X_t)}^4(\exp{(t+X_t)}^2-7\exp{(t+X_t)}+4)}{2(1+\exp{(t+X_t)})^7}.
\end{align*}
The corresponding numerical results are listed in Table \ref{cp:2}. Similar convergence results are shown as the above numerical tests.

\begin{table}[!htbp]
\caption{Numerical results of example 4.}\label{cp:2}
\begin{minipage}{\textwidth}\footnotesize%
\begin{center}%
\begin{tabular}{*{7}{c|}c}\hline
STEP& &N  &  $\abs{Y^0-Y_0}$   &  $\abs{Z^0-Z_0}$  &  $\abs{\Gamma^0-\Gamma_0}$  & $\abs{A^0-A_0}$  &  $T_{\mathrm rn}$ \\
\hline
\rb{-3em}{K=1} & \rb{-2.5em}{DT}
&   32 &  1.079E-03 &  2.877E-03 &  3.847E-03 &  3.427E-03 &
\\ \cline{3-7} &
&   64 &  5.374E-04 &  1.463E-03 &  2.008E-03 &  1.796E-03 &
\\ \cline{3-7} &
&  128 &  2.680E-04 &  7.389E-04 &  1.030E-03 &  9.216E-04 &
\\ \cline{3-7} &
&  256 &  1.338E-04 &  3.715E-04 &  5.225E-04 &  4.678E-04 &
\\ \cline{3-7} &
&  512 &  6.681E-05 &  1.864E-04 &  2.634E-04 &  2.359E-04 &
\\ \cline{2-7} & CR &&  1.00&  0.99&  0.97&  0.97& \rb{3em}{ 21.54s}
\footnote{\scriptsize -k 1 -ng 10 -ni 10 -e  F -sh  1.0 -mg -20.00 20.00 -$t_0$0.0 -$x_0$0.5 -T1.0}
\\ \hline
\rb{-3em}{K=2} & \rb{-2.5em}{DT}
&   32 &  5.131E-05 &  1.750E-05 &  2.281E-04 &  3.383E-04 &
\\ \cline{3-7} &
&   64 &  1.328E-05 &  4.871E-06 &  5.814E-05 &  8.806E-05 &
\\ \cline{3-7} &
&  128 &  3.388E-06 &  1.248E-06 &  1.477E-05 &  2.276E-05 &
\\ \cline{3-7} &
&  256 &  8.523E-07 &  3.171E-07 &  3.691E-06 &  5.698E-06 &
\\ \cline{3-7} &
&  512 &  2.139E-07 &  7.953E-08 &  9.282E-07 &  1.432E-06 &
\\ \cline{2-7} & CR &&  1.98&  1.95&  1.99&  1.97& \rb{3em}{ 70.10s}
\footnote{\scriptsize -k 2 -ng 10 -ni 10 -e  F -sh  1.0 -mg -20.00 20.00 -$t_0$0.0 -$x_0$0.5 -T1.0}
\\ \hline
\rb{-3em}{K=3} & \rb{-2.5em}{DT}
&   32 &  3.134E-06 &  6.034E-06 &  4.653E-06 &  2.167E-05 &
\\ \cline{3-7} &
&   64 &  4.058E-07 &  7.806E-07 &  5.391E-07 &  2.879E-06 &
\\ \cline{3-7} &
&  128 &  5.121E-08 &  9.992E-08 &  5.827E-08 &  3.605E-07 &
\\ \cline{3-7} &
&  256 &  6.407E-09 &  1.257E-08 &  7.059E-09 &  4.530E-08 &
\\ \cline{3-7} &
&  512 &  8.052E-10 &  1.574E-09 &  8.454E-10 &  5.398E-09 &
\\ \cline{2-7} & CR &&  2.98&  2.98&  3.11&  2.99& \rb{3em}{259.11s}
\footnote{\scriptsize -k 3 -ng 10 -ni 10 -e  F -sh  1.0 -mg -30.00 30.00 -$t_0$0.0 -$x_0$0.5 -T1.0}
\\ \hline
\end{tabular}\end{center}\end{minipage}
\end{table}

\subsection{Applications to stochastic optimal control}
We now show that one can solve stochastic optimal control problems in a 2FBSDEs way by using our multistep numerical schemes. To this end, let us consider the control problem, whose dynamic state equation is described by a forward SDE
\begin{equation}
\d X_t = b(t,X_t,\alpha_t) \d t + \sigma(t,X_t,\alpha_t) \d W_t,
\end{equation}
with the cost functional
\begin{equation}
J(\alpha) = \Ec{}{}{\int_0^T f(t,X_t,\alpha_t) \d t + g(X_t)}.
\end{equation}
The goal is to minimize the cost functional, i.e. Find $\alpha^\ast\in\U$ (where $\U$ contains all admissible controls, see e.g. \cite{MY} for the corresponding definition) satisfying
$$
V(t,x) :=J(t,x;\alpha^\ast) = \inf_{\alpha\in \U} J(t,x;\alpha).
$$
One can construct the corresponding HJB equation as
\begin{equation*}
\f{\p }{\p t}V(t,x) + \inf_{\alpha\in\U}\left\{\f{\sigma(t,x,\alpha)^2}2\f{\p^2 }{\p x^2}V(t,x) + b(t,x,\alpha)\f{ \p }{\p x}V(t,x)- f(t,x,\alpha)\right\} = 0.
\end{equation*}
and moreover, we have
\begin{equation}\label{ctrl:star}
\alpha^\ast(t,x) = \arg \inf_{\alpha\in\U}\left\{\f{\sigma(t,x,\alpha)^2}2\f{\p^2 V}{\p x^2} + b(t,x,\alpha)\f{ \p V}{\p x}- f(t,x,\alpha)\right\}.
\end{equation}
By inserting \eqref{ctrl:star} into the HJB equaion, one shows that the cost function satisfies
\begin{equation}
\f{\p }{\p t} V(t,x)+ G(t,x,\f{\p}{\p x} V(t,x),\f{\p^2}{\p x^2} V(t,x) ) = 0,
\end{equation}
where
$$
G(t,x,p,P) =  \f{\sigma(t,x,\alpha^\ast)^2}2 P + b(t,x,\alpha^\ast) p -  f(t,x,\alpha^\ast).
$$

For the above nonlinear PDEs, one can construct a corresponding 2FBSDEs. Then, one can solve the associate 2FBSDEs
using our numerical schemes to obtain $(X_t,Y_t,Z_t,\Gm_t),$ which yields the optimal state for the control problem, and finally, we can get obtain in view of (\ref{ctrl:star}) that
\begin{equation}
\alpha_t^\ast = g (X_t,Y_t,Z_t,\Gm_t)
\end{equation}
with $g(\cdot)$ being certain functions. This is a new approach dealing with the optimal control problem in a 2FBSDEs way. Now, we illustrate the idea by the following example.

\textbf{Tracking a particle under the microscope}:  consider the following system
\begin{equation}
\d X_t = \beta \alpha_t \d t + \sig \d W_t,
\end{equation}
where $X_t$ is the distance between the particle and the focus of the microscope,  $\beta \in \R$ is the gain in our servo loop and $\sig >0$  is the diffusion constant of the particle. We would like to keep the particle in focus, i.e.\ we expect that $X_t$ is as close to zero as possible. However, we have to introduce a power constraint on the control as well, as we cannot drive the servo motor with arbitrarily large input powers. We thus introduce the control cost
$$
J(\alpha) = \Ec{}{}{p\int_0^T X_t^2 \d t + q \int_0^T \alpha^2 \d t},
$$
where $p,q >0$ allows us to select the tradeoff between good tracking and low feedback power. One can construct the associate bellman equation
\begin{equation*}
\begin{aligned}
0 &= \f{\p }{\p t}V(t,x) + \inf_{\alpha\in \R}\left\{\f{\sig^2}2\f{\p^2}{\p x^2} V(t,x) + \beta\alpha\f{\p}{\p x}V(t,x) +p x^2 + q\alpha^2 \right\}\\
&=\f{\p }{\p t}V(t,x) + \f{\sig^2}2\f{\p^2 }{\p x^2}V(t,x) -\f{\beta}{4q}(\f{ \p }{\p x}V(t,x))^2 + px^2
\end{aligned}
\end{equation*}
with $V(T,x) =0$. Furthermore, one can show that the optimal control parameter (i.e., exact solution) is
\begin{equation}\label{ture_solution}
\alpha_t^\ast = -\f{\beta}{2q}\f{\p}{\p x} V(t,x).
\end{equation}
The classic numerical solution would relies on solving the above bellman equation.

We now solve the problem in a 2FBSDEs way, to this end, we first construct the corresponding 2FBSDEs as follows
\begin{eqsys}
\d X_t &= \beta c \d t + \sig \d W_t,\\
-\d Y_t &= (-\f{\beta^2}{4q\sig^2}Z_t^2 -\f{\beta c}{\sig} Z_t +pX_t^2) \d t - Z_t \d W_t,\\
\d Z_t &= A_t \d t + \Gm_t \d W_t.
\end{eqsys}
By solving this 2FBSDEs, we obtain the numerical solution in view of (\ref{ture_solution}) by
\begin{equation}
  \alpha^n  =\f{-\beta}{2q\sig}Z^n.
\end{equation}
In the numerical test, we set $\mu=0.1,\,r=0.03,\,\sig=0.5, \,c=0.1$. The numerical  results are shown in Table \ref{tab:ctrl1}. It can be seen from Table \ref{tab:ctrl1} that the approach is of high order accuracy, both for the 2FBSDE solution and the optimal control $\alpha$.

\begin{table}[!htbp]
\caption{Numerical results of stochastic optimal control.}\label{tab:ctrl1}
\begin{minipage}{\textwidth}\footnotesize%
\begin{center}%
\begin{tabular}{*{7}{c|}c}\hline
STEP& &N  &  $\abs{Y^0-Y_0}$   &  $\abs{Z^0-Z_0}$  &  $\abs{\Gamma^0-\Gamma_0}$  & $\abs{\alpha^0-\alpha^\ast}$  &  $T_{\mathrm rn}$ \\
\hline
\rb{-3em}{K=1} & \rb{-2.5em}{DT}
&   32 &  5.647E-03 &  5.856E-02 &  1.420E-02 &  1.952E-02 &
\\ \cline{3-7} &
&   64 &  2.788E-03 &  2.912E-02 &  7.075E-03 &  9.706E-03 &
\\ \cline{3-7} &
&  128 &  1.385E-03 &  1.452E-02 &  3.531E-03 &  4.839E-03 &
\\ \cline{3-7} &
&  256 &  6.903E-04 &  7.249E-03 &  1.764E-03 &  2.416E-03 &
\\ \cline{3-7} &
&  512 &  3.446E-04 &  3.622E-03 &  8.814E-04 &  1.207E-03 &
\\ \cline{2-7} & CR &&  1.01&  1.00&  1.00&  1.00& \rb{3em}{  2.05s}
\footnote{\scriptsize -k 1 -ng 10 -ni 8 -e  F -sh  1.0 -mg -20.00 20.00 -$t_0$0.0 -$x_0$5.0 -T1.0}
\\ \hline
\rb{-3em}{K=2} & \rb{-2.5em}{DT}
&   32 &  1.438E-03 &  8.777E-04 &  5.823E-05 &  2.926E-04 &
\\ \cline{3-7} &
&   64 &  3.660E-04 &  2.192E-04 &  1.471E-05 &  7.308E-05 &
\\ \cline{3-7} &
&  128 &  9.230E-05 &  5.478E-05 &  3.696E-06 &  1.826E-05 &
\\ \cline{3-7} &
&  256 &  2.317E-05 &  1.369E-05 &  9.264E-07 &  4.564E-06 &
\\ \cline{3-7} &
&  512 &  5.806E-06 &  3.422E-06 &  2.319E-07 &  1.141E-06 &
\\ \cline{2-7} & CR &&  1.99&  2.00&  1.99&  2.00& \rb{3em}{  6.74s}
\footnote{\scriptsize -k 2 -ng 10 -ni 8 -e  F -sh  1.0 -mg -20.00 20.00 -$t_0$0.0 -$x_0$5.0 -T1.0}
\\ \hline
\rb{-3em}{K=3} & \rb{-2.5em}{DT}
&   32 &  3.439E-06 &  4.495E-06 &  1.848E-06 &  1.498E-06 &
\\ \cline{3-7} &
&   64 &  4.758E-07 &  5.315E-07 &  2.276E-07 &  1.772E-07 &
\\ \cline{3-7} &
&  128 &  6.243E-08 &  6.456E-08 &  2.824E-08 &  2.152E-08 &
\\ \cline{3-7} &
&  256 &  7.992E-09 &  7.953E-09 &  3.517E-09 &  2.651E-09 &
\\ \cline{3-7} &
&  512 &  1.011E-09 &  9.869E-10 &  4.387E-10 &  3.290E-10 &
\\ \cline{2-7} & CR &&  2.94&  3.04&  3.01&  3.04& \rb{3em}{ 19.44s}
\footnote{\scriptsize -k 3 -ng 10 -ni 8 -e  F -sh  1.0 -mg -20.00 20.00 -$t_0$0.0 -$x_0$5.0 -T1.0}
\\ \hline
\rb{-3em}{K=4} & \rb{-2.5em}{DT}
&   32 &  4.583E-07 &  3.240E-07 &  3.570E-08 &  1.080E-07 &
\\ \cline{3-7} &
&   64 &  2.986E-08 &  2.005E-08 &  2.264E-09 &  6.682E-09 &
\\ \cline{3-7} &
&  128 &  1.902E-09 &  1.246E-09 &  1.425E-10 &  4.154E-10 &
\\ \cline{3-7} &
&  256 &  1.200E-10 &  7.767E-11 &  8.933E-12 &  2.589E-11 &
\\ \cline{3-7} &
&  512 &  7.535E-12 &  4.847E-12 &  5.592E-13 &  1.616E-12 &
\\ \cline{2-7} & CR &&  3.97&  4.01&  3.99&  4.01& \rb{3em}{ 57.62s}
\footnote{\scriptsize -k 4 -ng 10 -ni 8 -e  F -sh  1.0 -mg -20.00 20.00 -$t_0$0.0 -$x_0$5.0 -T1.0}
\\ \hline
\rb{-3em}{K=5} & \rb{-2.5em}{DT}
&   32 &  4.302E-09 &  2.299E-10 &  8.347E-10 &  7.662E-11 &
\\ \cline{3-7} &
&   64 &  1.512E-10 &  1.904E-11 &  2.438E-11 &  6.346E-12 &
\\ \cline{3-7} &
&  128 &  4.996E-12 &  7.749E-13 &  7.352E-13 &  2.583E-13 &
\\ \cline{3-7} &
&  256 &  1.604E-13 &  2.698E-14 &  2.256E-14 &  8.995E-15 &
\\ \cline{3-7} &
&  512 &  5.080E-15 &  8.862E-16 &  6.985E-16 &  2.954E-16 &
\\ \cline{2-7} & CR &&  4.93&  4.54&  5.05&  4.54& \rb{3em}{ 94.61s}
\footnote{\scriptsize -k 5 -ng 10 -ni 12 -e  F -sh  1.0 -mg -30.00 30.00 -$t_0$0.0 -$x_0$5.0 -T1.0}
\\ \hline
\end{tabular}\end{center}\end{minipage}
\end{table}

\section{Conclusions}\label{cons}
We have extended our multistep schemes in \cite{ZFZ2014} to the multistep schemes 3 and 5 for solving the second order FBSDEs. The key feature of the proposed multistep schemes is that the Euler method is used to discrete the forward SDE, which dramatically reduces the entire computational complexity. Furthermore, it is shown that the quantities of interest (e.g., the solution tuple $(Y_t, Z_t, A_t, \Gamma_t)$ in the 2FBSDEs) are still of high order accuracy. Several numerical examples are presented to show the effective of the proposed numerical schemes.  Applications of our numerical schemes for stochastic optimal control problems are also discussed.

There are, however, some other related topics that need to be investigated:
\begin{itemize}
\item High dimensional problems. Note the methods here can be easily extended to high dimensional problems. However, we have proposed the local Lagrange interpolation methods here in our schemes. For high dimensional problems, this would results in the tensor Lagrange interpolation methods, which may be time consuming. Thus, we would suggest more feasible techniques such as the sparse grid interpolation, RBF interpolation, etc. This would be part of our future studies.

\item Rigorous stability and convergence analysis. This is also our ongoing project.
\end{itemize}

\end{document}